\documentclass[a4paper,12pt]{article}

\usepackage{amsmath}
\usepackage{amsfonts}
\usepackage{a4wide}

\newenvironment{proof}{\textbf{Proof:}} {\hspace{\stretch{1}}$\Box$}

\newcommand{\fbM}[4]{\left(\begin{array}{c|c} 
#1 & #2\\
\hline
#3 & #4
\end{array}\right)}
\newcommand{\fcv}[4]{\left(\begin{array}{c c|c c} 
#1 & #2 & #3 & #4
\end{array}\right)}
\newcommand{\frv}[4]{\left(\begin{array}{c} 
#1\\ #2\\
\hline
#3\\ #4\\
\end{array}\right)}
\newcommand{\subM}[4]{\begin{array}{c c} 
#1 & #2\\
#3 & #4
\end{array}}

\newcommand{\norm}[3]{\left\|#1\right\|^{#2}_{#3}}  

\newcommand{\Pmp}[1]{P^\mp(#1)} 
\newcommand{\Ppm}[1]{P^\pm(#1)} 

\newcommand{\diag}[2]{\textrm{diag}(#1,#2)} 

\newcommand{\tmu}{\tilde\mu}
\newcommand{\tga}{\tilde\gamma}
\newcommand{\hmu}{\hat\mu}
\newcommand{\hga}{\hat\gamma}
\newcommand{\hz}{\hat\zeta}
\newcommand{\chz}{\check\zeta}
\makeatletter
\def\@begintheorem#1#2{\trivlist \item[\hskip \labelsep{\sc #1\ #2.}]\it}
\def\@opargbegintheorem#1#2#3{\trivlist
      \item[\hskip \labelsep{\sc #1\ #2.\ (#3)}]\it}
\def\@sect#1#2#3#4#5#6[#7]#8{\ifnum #2>\c@secnumdepth
     \let\@svsec\@empty\else
     \refstepcounter{#1}\edef\@svsec{\csname the#1\endcsname.\hskip 1em}\fi
     \@tempskipa #5\relax
      \ifdim \@tempskipa>\z@
        \begingroup #6\relax
          \@hangfrom{\hskip #3\relax\@svsec}{\interlinepenalty \@M #8\par}%
        \endgroup
       \csname #1mark\endcsname{#7}\addcontentsline
         {toc}{#1}{\ifnum #2>\c@secnumdepth \else
                      \protect\numberline{\csname the#1\endcsname}\fi
                    #7}\else
        \def\@svsechd{#6\hskip #3\relax  
                   \@svsec #8\csname #1mark\endcsname
                      {#7}\addcontentsline
                           {toc}{#1}{\ifnum #2>\c@secnumdepth \else
                             \protect\numberline{\csname the#1\endcsname}\fi
                       #7}}\fi
     \@xsect{#5}}
\def\section{\@startsection {section}{1}{\z@}{-1.5ex plus-1ex minus
    -.2ex}{-2.5ex plus.2ex}{\reset@font\bf}}
\def\subsection{\@startsection{subsection}{2}{\z@}{-3.25ex plus-1ex
    minus-.2ex}{-1.5ex plus.2ex}{\reset@font\sl}}

\makeatother

\newtheorem{predf}{Definition}[section]

\newtheorem{prex}[subsection]{Example}

\newtheorem{prop}[subsection]{Proposition}
\newtheorem{lemma}[subsection]{Lemma}

\newcommand{\R}{{\Bbb R}}

\def \beq {\begin {eqnarray}}
\def \eeq {\end {eqnarray}}
\def \ba {\begin {eqnarray*}}
\def \ea {\end  {eqnarray*}}
\def \p {\partial}
\def \tilde {\widetilde}
\def \hat {\widehat}
     \newcommand{\be}{\begin{equation}}   
\newcommand{\ee}{\end{equation}}

\newcommand{\bea}{\begin{eqnarray}}  
\newcommand{\eea}{\end{eqnarray}}
\newcommand{\bean}{\begin{eqnarray*}}
\newcommand{\eean}{\end{eqnarray*}}
\newcommand{\abs}[1]{\lvert #1 \rvert}          

\newcommand{\closure}[1]{\overline{#1}}

\begin{document}
\title{Inverse boundary value problem for  Maxwell equations with local data
\footnote{Supported by the Academy of Finland under
CoE--project 213476}
\footnote{The first author was also supported by Ministerio de Ciencia e Innovaci\'on de Espa\~na, MTM2005-07652-C02-01}}

\author{
Pedro Caro\footnote{Department of Mathematics, Universidad Aut\'onoma de Madrid, 28049 Madrid, Spain}, Petri Ola\footnote{Department of Mathematics and Statistics, P.O. Box 68, 00014 University of Helsinki, Finland} and 
Mikko Salo\footnote{Department of Mathematics and Statistics, P.O. Box 68, 00014 University of Helsinki, Finland}
\\
}
\date{February 6, 2009}

\maketitle

\begin{abstract}
\noindent We prove a uniqueness theorem for an inverse boundary value problem for the Maxwell system with boundary data assumed known only in part of the boundary. We assume that the inaccessible part of the boundary is either part of a plane, or part of a sphere. 
This work generalizes the results obtained by Isakov \cite{I} for the Schr\"odinger equation to Maxwell equations.

\end{abstract}

\section*{Introduction.} Let \(\Omega\subset\R^3\) be a bounded domain with $C^{1,1}$ boundary, and let $\varepsilon$, $\mu$, $\sigma$ be $C^2$ functions in $\closure{\Omega}$ ($\varepsilon$ is the permittivity, $\mu$ the permeability, and $\sigma$ the conductivity). We will assume that the coefficients satisfy the positivity conditions 
\begin{equation} \label{PositivityConditions}
\gamma = \varepsilon + i\sigma/\omega, \ \varepsilon>0, \, \mu >0, \, \sigma \geq 0 \ {\rm in} \, \overline\Omega.
\end{equation}
Let \(D= -i\nabla\), let \(\nu\) be the exterior unit normal to \(\p\Omega\), and consider the time-harmonic Maxwell equations for the electric field $E$ and magnetic field $H$ in \(\Omega\),
\begin{equation}\label{Maxwell_system_in_Omega_Intro}
\left\{
\begin{array}{l}
D\wedge H+\omega\gamma E=0,\\
D\wedge E-\omega\mu H=0,
\end{array}
\right.
\end{equation}
with the boundary condition
\be
\nu \wedge H = a \ {\rm on}\ \p\Omega.
\ee
This is the {\em magnetic boundary value problem} for the Maxwell equations.
Here we use \('\wedge '\) to denote the vector product  in \(\R ^3\), and \(\nabla \wedge F\) is the curl of the vector field \(F\). When posed in correct function spaces this problem admits a unique solution \( (E,H)\) when the angular frequency \(\omega >0\) does not belong to a discrete set of resonant frequencies. The {\em impedance map} \(\Lambda (\gamma, \mu)\) is then defined by
\[
\Lambda (\gamma, \mu): \  a\mapsto \nu \wedge E|_{\partial\Omega}.
 \]
In analogue with the Calder\'on problem of electrostatics (see \cite{C}, \cite{N} and \cite{S-U}), we consider the inverse problem of determining the electromagnetic parameters \(\gamma\) and \(\mu\) from the boundary measurement \(\Lambda (\gamma, \mu)\). Especially, in this work we assume that we can measure the values of \(\Lambda (\gamma, \mu)(a)\) only on a nonempty open subset \(\Gamma\) of \(\p\Omega\), and only for tangential boundary fields \(a\) supported in \(\Gamma\). We further assume that the inaccessible part of the boundary \(\Gamma _0 = \overline{\p\Omega\setminus \Gamma}\) is either part of a sphere, or part of a plane. Similar results were proved by Isakov \cite{I} for the inverse boundary value problem for the Schr\"odinger operator. Our work extends this method to Maxwell equations. 
\smallskip

Let us now formulate the main results of this paper in a precise way. Given \(\Omega\) as above, we define the spaces 
\begin{gather*}
H^1_{\rm Div}(\Omega) = \{ F\in H^1(\Omega)^3; \, {\rm Div}(\nu \wedge F) \in H^{1/2}(\partial\Omega)\}, \\
TH_{\rm Div}^{1/2}(\p\Omega) = \{ a\in H^{1/2}(\p\Omega)^3; \, \nu\cdot a = 0, \, {\rm Div}(a)\in H^{1/2}(\p\Omega)\},
\end{gather*}
where \({\rm Div}\) denotes the surface divergence on \(\p\Omega\) (see for example \cite{C-P} and \cite{O-P-S1} for more information). Provided that $\mu$ and $\gamma$ are coefficients in $C^2(\closure{\Omega})$ which satisfy \eqref{PositivityConditions}, then for \(\omega\) outside a discrete set of resonant frequencies, the magnetic boundary value problem \eqref{Maxwell_system_in_Omega_Intro} has a unique solution in $H^1_{\rm Div}(\Omega) \times H^1_{\rm Div}(\Omega)$ for any $a \in TH_{\rm Div}^{1/2}(\p\Omega)$ (see \cite{O-P-S1}). The impedance map is well defined, and 
\[
\Lambda (\gamma, \mu): TH_{\rm Div}^{1/2}(\p\Omega) \to TH_{\rm Div}^{1/2}(\p\Omega).
\]
Our first main result is the following.
\medskip

{\sc Theorem 1.}\, {\it
Assume that \(\Omega\subset\{\R ^3; \, x_3 <0\}\) is a $C^{1,1}$ domain. Let
\(\Gamma _0 = \p\Omega\cap \{ x_3 = 0\}\) and  \(\Gamma = \overline{\p\Omega\setminus \Gamma _0}.\)
Assume we have two pairs of electromagnetic parameters \((\gamma_j, \mu_j) \in C^4 (\overline \Omega)\times
C^4 (\overline \Omega)\), \(j= 1,\, 2\),  satisfying the positivity conditions \eqref{PositivityConditions} and the following boundary conditions: 
\begin{gather}
\gamma _1 = \gamma _2, \ \mu _1 = \mu _2 \text{ up to order one on \(\Gamma\), and } \label{BoundaryGammaAssumption} \\ 
\left\{ \begin{array}{l} \text{there exist $C^4$ extensions of $\gamma_j$ and $\mu_j$ into $\R^3$ which are } \\[1pt] \text{invariant under reflection across the plane $\{x_3 = 0\}$.} \end{array} \right. \label{BoundaryExtensionAssumption}
\end{gather}
Suppose that $\omega$ is not a resonant frequency for $(\gamma_1,\mu_1)$ or $(\gamma_2,\mu_2)$. If the impedance maps restricted to \(\Gamma\) coincide,
\[
\Lambda(\gamma_1, \mu_1)(a)|_{\Gamma} = \Lambda(\gamma_2, \mu_2)(a)|_{\Gamma} \ \text{for all $a \in TH^{1/2}_{\rm Div}(\partial \Omega)$ with $\text{supp}(a) \subset \Gamma$},
\]
then the electromagnetic parameters are equal, i.e. \(\gamma _1  = \gamma _2\) and \(\mu _1 = \mu _2\) in \(\Omega\).}
\medskip

Before formulating the second main theorem of our paper, we want to comment on the assumptions in Theorem 1. The condition \eqref{BoundaryGammaAssumption} is not important, and we expect that a suitable boundary determination result would allow to remove this condition completely (in the case of $C^{\infty}$ coefficients such a result is given in \cite{J-M}). The condition \eqref{BoundaryExtensionAssumption} comes from the method of proof, which requires that the Maxwell equations are invariant under reflection across $\{x_3 = 0\}$. The following is a necessary condition for \eqref{BoundaryExtensionAssumption} to hold:
\begin{equation} \label{BoundaryExtensionNecessary}
\left(\frac{\partial}{\partial \nu}\right)^l \gamma_1 = \left(\frac{\partial}{\partial \nu}\right)^l \gamma_2 = 0, \ \left(\frac{\partial}{\partial \nu}\right)^l \mu_1 = \left(\frac{\partial}{\partial \nu}\right)^l \mu_2 = 0 \text{ on $\Gamma_0$ for $l = 1, 3$}.
\end{equation}

\smallskip

For the sake of definiteness, we have stated Theorem 1 in terms of the impedance map. However, the proof of Theorem 1 extends to the case of restricted Cauchy data sets 
\begin{equation*}
C_{\Gamma}(\gamma,\mu) = \{ (\nu \wedge E|_{\Gamma}, \nu \wedge H|_{\Gamma}) ; (E,H) \in H^2(\Omega)^3 \times H^2(\Omega)^3 \text{ solves } \eqref{Maxwell_system_in_Omega_Intro} \}.
\end{equation*}
In this setup $\Omega$ can be a Lipschitz domain and $\omega > 0$ any number, and the assumption 
\begin{equation*}
C_{\Gamma}(\gamma_1,\mu_1) = C_{\Gamma}(\gamma_2,\mu_2)
\end{equation*}
for two pairs $(\gamma_j,\mu_j)$ satisfying the conditions in Theorem 1 will imply that $\gamma_1 = \gamma_2$ and $\mu_1 = \mu_2$ in $\Omega$. This result allows domains with transversal intersections
of \(\p\Omega\) and \(\{x_3 =0\}\), in which case \eqref{BoundaryExtensionNecessary} is also sufficient for \eqref{BoundaryExtensionAssumption}. It is well known (see for example \cite{M-M}, \cite{S1} and \cite{S2} and the references given in these articles) that Maxwell equations may not admit \(H^1\)-solutions even with boundary data in \(H^{1/2}(\p\Omega)\) if the domain is non convex: the solution may have a finite dimensional singular part near the conical singularity, so the image of the impedance map might be slightly larger than just \(TH_{\rm Div}^{1/2}(\p\Omega)\). However, all solutions needed in the proof are restrictions to \(\Omega\) of \(H^2\)-solutions defined in a neighborhood of \(\Omega\) (actually they are defined in the  whole \(\R ^3\)), and hence it is enough to restrict only to the part of the graph of \(\Lambda (\gamma, \mu)\) that is actually needed.
\smallskip

Let us now formulate our second main theorem:
\medskip

{\sc Theorem 2.}\,{\it 
Assume \(\Omega\subset B_0\) is a $C^{1,1}$ domain included in an open ball \( B_0\) of \(\R ^3\) of positive radius,  and let
\(\Gamma _0 = \p\Omega\cap \p B_0,\ \ \Gamma_0 \neq \partial B_0, \ \ \Gamma = \overline{\p\Omega\setminus \Gamma _0}.\)
Assume we have two pairs of electromagnetic parameters \((\gamma _j , \mu _j) \in C^4 (\overline \Omega)\times
C^4 (\overline \Omega)\), \(j= 1,\, 2\),  satisfying the conditions of Theorem 1 but with respect to \(\p B_0\)
instead of \(\{x_3 = 0\}\).
If $\omega$ is not a resonant frequency and if the impedance maps restricted to \(\Gamma\) coincide,
\[
\Lambda(\gamma_1, \mu_1)(a)|_{\Gamma} = \Lambda(\gamma_2, \mu_2)(a)|_{\Gamma} \ \text{for all $a \in TH^{1/2}_{\rm Div}(\partial \Omega)$ with $\text{supp}(a) \subset \Gamma$},
\]
then the electromagnetic parameters are equal, i.e. \(\gamma _1  = \gamma _2\) and \(\mu _1 = \mu _2\) in \(\Omega\).}
\medskip

The unique recovery of electromagnetic parameters from the scattering amplitude was first proven in \cite{C-P} under the assumption that the magnetic permeability \(\mu\) is a constant. The unique recovery of general \(C^2\)-parameters \(\gamma\) and \(\mu\) from full boundary data was then proved in \cite{O-P-S1}, and later simplified in \cite{O-S}, and the uniqueness for the inverse scattering problem  was proved in \cite{Sa}. Boundary determination results were given in \cite{McD1} and \cite{J-M}, and the more general chiral media was studied in \cite{McD2}.  For a slightly more general approach and more background information, see also the review article \cite{O-P-S2}. Note that in this work we need, for technical reasons, two more derivatives of the parameters compared to the full data problems.
\smallskip

Inverse problems with partial data for scalar elliptic equations have attracted considerable attention recently. There are two main approaches: the first uses Carleman estimates \cite{B-U}, \cite{K-S-U} and the second uses reflection arguments \cite{I}. In the first approach $\Omega$ can be any domain but one needs to measure part of Cauchy data in a small set $\Gamma \subset \partial \Omega$ and the other part in a neighborhood of $\partial \Omega \setminus \Gamma$. In the second approach it is enough to measure Cauchy data on a fixed small set $\Gamma$, but one has the restriction that $\partial \Omega \setminus \Gamma$ has to be part of a plane or sphere.
\smallskip

In this paper, we apply Isakov's reflection method \cite{I} to the Maxwell system to prove Theorems 1 and 2. As far as we know, these are the first partial data results for an elliptic system which are analogous to \cite{I} (an analog of \cite{K-S-U} for a Dirac system was recently proved in \cite{SaT}). The strategy of the proof is similar to \cite{O-P-S1} and \cite{O-S}: we construct special solutions to the Maxwell system, insert them in a suitable integral identity, and recover the coefficients from an asymptotic limit.
\smallskip

However, several issues arise when trying to combine this strategy with the reflection method. The integral identities in \cite{O-P-S1} and \cite{O-S} do not seem to go well together with partial data, and the final identity which we use for recovering the coefficients is new. In this identity, many terms survive in the asymptotic limit, and one needs to carefully manipulate these terms to determine the coefficients. Also, there are cross terms involving products of the original and reflected solutions. In \cite{I} these cross terms were handled by the Riemann-Lebesgue lemma. In our case we need to prove a substitute for the Riemann-Lebesgue lemma, and here the $C^4$-regularity of the coefficients is used.
\smallskip

The structure of the paper is as follows. In the first section we augment the Maxwell equations into a Dirac type elliptic system and prove a crucial factorization result for a second order matrix Schr\"odinger operator. This follows in principle \cite{O-S}, but with somewhat different  notations and some simplifications. In order to keep the article self-contained, we have included complete proofs. In the second section we construct the complex geometric optics solutions  (CGO solutions for short) needed in the proof. In this section, following Isakov \cite{I}, we also perform the reflection across \(\{x_3 = 0\}\) and analyze the behavior of CGO solutions for Maxwell equations under this operation. 
In the third section we  prove the integration by parts formula needed in our argument, and in the  fourth section we perform the asymptotical computations needed in our proof. Theorem 1 is then proved in section five, and finally using the Kelvin transform, the proof of Theorem 2 is reduced to Theorem 1 in the last section. Note that with the exception of the last section, we use the vector field formulation for our equations. Only when analyzing the behavior of Maxwell equations under the Kelvin transform does it become useful to revert to a formulation using differential forms. However, we feel that the major part of this work is easier to read when written in vector field notation. 

\textbf{Acknowledgement:} The first author would like to thank Lassi P\"aiv\"arinta and the Inverse Problems group at the University of Helsinki for their hospitality.

\section{Augmented Maxwell system and reduction to Schr\"odinger equation.}\label{sec:augmented_Max}
In this section we recall an argument used in \cite{O-S}. Roughly, the idea is to augment Maxwell equations in such a way that, with a suitable rescaling of the fields, the system can be transformed into a Schr\"odinger equation. This argument can be written in a slightly more general setting than the one described above, see \cite{O-P-S2}. Let $U\subset\mathbb{R}^3$ be an an open set along this section. Consider the time-harmonic Maxwell system in $U$ given by
\begin{equation}\label{Maxwell_system_in_O}
\left\{
\begin{array}{l}
D\wedge H+\omega\gamma E=0\\
-D\wedge E+\omega\mu H=0
\end{array}
\right.
\end{equation}
where $\omega>0$ is a fixed angular frequency, $\gamma=\varepsilon+i\sigma/\omega$ and $\mu,\varepsilon,\sigma\in C^2(\overline{U})$ satisfy 
\[
\mu, \, \varepsilon >0, \ \ \sigma \geq 0 \ \mbox{in \(\overline U\)}.
\]
Taking the divergence of both equations in \eqref{Maxwell_system_in_O} we obtain the additional equations
\begin{equation*}
\left\{
\begin{array}{l}
D\cdot(\gamma E)=0\\
D\cdot(\mu H)=0.
\end{array}
\right.
\end{equation*}
These will be useful in writing the Maxwell equations as an elliptic system. Namely, denoting $\alpha=\textrm{log }\gamma$ and $\beta=\textrm{log }\mu$ we may combine these into four equations
\begin{equation*}
\left\{
\begin{array}{l}
D\cdot E+D\alpha\cdot E=0\\
-D\wedge E+\omega\mu H=0\\
D\cdot H+D\beta\cdot H=0\\
D\wedge H+\omega\gamma E=0.
\end{array}
\right.
\end{equation*}
This gives rise to the $8\times 8$ system
\begin{equation*}
\left[
\left(\begin{array}{c c c c}
* & 0 & * & D\cdot\\
* & 0 & * & -D\wedge\\
* & D\cdot & * & 0\\
* & D\wedge & * & 0\\
\end{array}\right)
+\left(\begin{array}{c c c c}
* & 0 & * & D\alpha\cdot\\
* & \omega\mu I_3 & * & 0\\
* & D\beta\cdot & * & 0\\
* & 0 & * & \omega\gamma I_3\\
\end{array}\right)
\right]\left(\begin{array}{c}
0\\ H\\ 0\\ E\\
\end{array}\right)=0.
\end{equation*}
Here $*$ means that in these positions we may insert anything and would still get the same equations.\\

We wish to consider $8$-vectors $X^t=\fcv{\Phi}{H^t}{\Psi}{E^t}$, where $\Psi$ and $\Phi$ are additional scalar fields. We also want the system to be elliptic, in fact, we want the principal part to be a Dirac-type operator. To this end, consider two $4\times 4$ operators (acting on $4$-vectors $\left( \begin{array}{cc} \Phi & H^t \end{array} \right)^t$ for instance) 
\begin{equation*}
P_+(D)=
\left(\begin{array}{c c}
0 & D\cdot\\
D & D\wedge
\end{array}\right),\quad
P_-(D)=
\left(\begin{array}{c c}
0 & D\cdot\\
D & -D\wedge
\end{array}\right).
\end{equation*}
These operators satisfy
\begin{gather*}
P_+(D)P_-(D)=P_-(D)P_+(D)=-\Delta I_4,\\
P_+(D)^*=P_-(D),\quad P_-(D)^*=P_+(D).
\end{gather*}
Motivated by the Pauli-Dirac operator in $\mathbb{R}^3$, we choose the principal part to have the block form
\begin{equation*}
\Pmp{D}=\left(\begin{array}{c c}
0 & P_-(D)\\
P_+(D) & 0
\end{array}\right).
\end{equation*}
Observe that $\Pmp{D}\Pmp{D}=-\Delta I_8$ and $\Pmp{D}^*=\Pmp{D}$.
The following notations will also be used in the rest of  the work. Let us denote
\begin{equation}
\Pmp{a,b}=\fbM{0}{P_-(b)}{P_+(a)}{0}, \qquad \Ppm{a,b}=\fbM{0}{P_+(b)}{P_-(a)}{0},
\end{equation}
where $a,b\in\mathbb{C}^3$; and
\[\diag{A}{B}=\fbM{A}{0}{0}{B},\]
where $A$ and $B$ are $4\times 4$ diagonal $\mathbb{C}$-matrices. 
When $A=\lambda_1I_4$ and $B=\lambda_2I_4$ with $\lambda_1,\lambda_2\in\mathbb{C}$ let us write
\[\diag{\lambda_1}{\lambda_2}=\fbM{\lambda_1I_4}{0}{0}{\lambda_2I_4}.\]
Note that $\Pmp{a,b}\Pmp{b,a}=\diag{b\cdot b}{a\cdot a}$ and $\Ppm{a,b}\Ppm{b,a}=\diag{b\cdot b}{a\cdot a}$ for any $a,b\in\mathbb{C}^3$.
We shall write for short $\Pmp{a}$ and $\Ppm{a}$ when $b=a$, hence $\Pmp{a}\Pmp{a}=(a\cdot a)I_8$ and $\Ppm{a}\Ppm{a}=(a\cdot a)I_8$ for any $a\in\mathbb{C}^3$. From these last equalities follow the anticommutation formulas
\begin{gather}\label{pseudo-commutation1}
\Pmp{a}\Pmp{b}=-\Pmp{b}\Pmp{a}+2(a\cdot b)I_8\\
\label{pseudo-commutation2}
\Ppm{a}\Ppm{b}=-\Ppm{b}\Ppm{a}+2(a\cdot b)I_8,
\end{gather}
for $a,b\in\mathbb{C}^3$. Moreover, the commutation formula
\begin{equation}\label{commutation}
\Ppm{a}\Pmp{b}=\Pmp{b}\Ppm{a}
\end{equation}
holds for $a,b\in\mathbb{C}^3$. This formula follows easily from $P_+(a)P_+(b)=(P_+(a)P_+(b))^t=P_-(b)P_-(a)$. More generally
\begin{equation}\label{commutation_gen}
\Ppm{a,b}\Pmp{c}=\Pmp{c}\Ppm{b,a},
\end{equation}
for $a,b,c\in\mathbb{C}^3$.

For the potential there is seemingly considerable freedom in the choice of the entries marked with $*$. However, we later wish to reduce the Maxwell equations into a matrix Schr\"odinger equation with no first order term so a proper choice of the extra entries will be crucial. We will take
\[V_{\mu,\gamma}=\left(\begin{array}{c c c c}
\omega\mu & 0 & 0 & D\alpha\cdot\\
0 & \omega\mu I_3 & D\alpha & 0\\
0 & D\beta\cdot & \omega\gamma & 0\\
D\beta & 0 & 0 & \omega\gamma I_3\\
\end{array}\right).\]
The augmented Maxwell system is then
\[
(\Pmp{D}+V_{\mu,\gamma})X=0 \ \ {\rm in}Ê\   U.
\]
The solutions of this system for which $\Phi=\Psi=0$ correspond to the solutions of the original Maxwell system in $U$.

We may write the potential \(V_{\mu,\gamma}\) above  in the form
\begin{equation}\label{potential_V}
V_{\mu,\gamma}=\omega\diag{\mu}{\gamma}+\frac{1}{2} (\Ppm{D\beta,D\alpha}+\Pmp{D\beta,D\alpha}).
\end{equation}
Next we rescale the \(X\) as follows. Let 
\[
Y^t=\fcv{Y^\Phi}{(Y^H)^t}{Y^\Psi}{(Y^E)^t}
\]
be defined by
\[
Y=\diag{\mu^{1/2}}{\gamma^{1/2}}X.
\]
Then 
\[(\Pmp{D}+V_{\mu,\gamma})X=\diag{\gamma^{-1/2}}{\mu^{-1/2}}(\Pmp{D}+W_{\mu,\gamma})Y\]
where
\[W_{\mu,\gamma}=
\fbM{\kappa I_4}{\frac{1}{2}P_+(D\alpha)}{\frac{1}{2}P_-(D\beta)}{\kappa I_4}
=\kappa I_8+\frac{1}{2}\Ppm{D\beta,D\alpha}\]
and $\kappa=\omega(\gamma\mu)^{1/2}$. Here we used that
\[\Pmp{D}(\diag{\mu^{-1/2}}{\gamma^{-1/2}}Y)= \diag{\gamma^{-1/2}}{\mu^{-1/2}}[\Pmp{D}-\frac{1}{2}\Pmp{D\beta,D\alpha}]Y.\]
Hence
\[(\Pmp{D}+V_{\mu,\gamma})X=0 \quad \Leftrightarrow \quad (\Pmp{D}+W_{\mu,\gamma})Y=0\]
when $X$ and $Y$ are related as above. This scaling is motivated by the following result.\\

\begin{lemma}\label{lemma:schroedinger}
One has
\begin{gather*}
(\Pmp{D}+W_{\mu,\gamma})(\Pmp{D}-W^t_{\mu,\gamma})=-\Delta I_8+\tilde{Q}_{\mu,\gamma}\\
(\Pmp{D}-W^t_{\mu,\gamma})(\Pmp{D}+W_{\mu,\gamma})=-\Delta I_8+\tilde{Q}'_{\mu,\gamma}
\end{gather*}
where the matrix potentials are given by
\setlength\arraycolsep{2pt}
\begin{eqnarray*}
\tilde{Q}_{\mu,\gamma}&=&\frac{1}{2}
\fbM{\subM{\Delta\alpha}{0}{0}{2\nabla^2\alpha-\Delta\alpha I_3}}{0}
{0}{\subM{\Delta\beta}{0}{0}{2\nabla^2\beta-\Delta\beta I_3}}\\
& &-\fbM{(\kappa^2+\frac{1}{4}(D\alpha)^2)I_4}{\subM{0}{2D\kappa\cdot}{2D\kappa}{0}}
{\subM{0}{2D\kappa\cdot}{2D\kappa}{0}}{(\kappa^2+\frac{1}{4}(D\beta)^2)I_4}
\end{eqnarray*}
and
\setlength\arraycolsep{2pt}
\begin{eqnarray*}
\tilde{Q}'_{\mu,\gamma}&=&-\frac{1}{2}
\fbM{\subM{\Delta\beta}{0}{0}{2\nabla^2\beta-\Delta\beta I_3}}{0}
{0}{\subM{\Delta\alpha}{0}{0}{2\nabla^2\alpha-\Delta\alpha I_3}}\\
& &-\fbM{(\kappa^2+\frac{1}{4}(D\beta)^2)I_4}{\subM{0}{0}{0}{2D\kappa\wedge}}
{\subM{0}{0}{0}{-2D\kappa\wedge}}{(\kappa^2+\frac{1}{4}(D\alpha)^2)I_4}
\end{eqnarray*}
with $\nabla^2f=(\partial^2_{x_j,x_k}f)^3_{j,k=1}$.
\end{lemma}

\begin{proof}
We have
\[(\Pmp{D}+W_{\mu,\gamma})(\Pmp{D}-W_{\mu,\gamma}^t)=-\Delta I_8+W_{\mu,\gamma}\Pmp{D} -\Pmp{D}W_{\mu,\gamma}^t-W_{\mu,\gamma}W_{\mu,\gamma}^t.\]
We claim that
\begin{gather*}
W_{\mu,\gamma}\Pmp{D}-\Pmp{D}W_{\mu,\gamma}^t=\\
=\frac{1}{2}\fbM{\subM{\Delta\alpha}{0}{0}{2\nabla^2\alpha-\Delta\alpha I_3}}{0}
{0}{\subM{\Delta\beta}{0}{0}{2\nabla^2\beta-\Delta\beta I_3}}-\Pmp{D\kappa}.\\
\end{gather*}
This can be proved by a direct calculation. However, we will use symbol calculus which gives a slightly more elegant proof. In standard (left) quantization, since the symbol of $\Pmp{D}$ is $\Pmp{\xi}$ and since $\partial^2_{\xi_j}\Pmp{\xi}=\partial_{\xi_j} W_{\mu,\gamma}=0$ for $j=1,2,3$, we have
\setlength\arraycolsep{2pt}
\begin{eqnarray*}
W_{\mu,\gamma}\Pmp{D}-\Pmp{D}W_{\mu,\gamma}^t&=&
\textrm{Op}(W_{\mu,\gamma}\Pmp{\xi}-\Pmp{\xi}W_{\mu,\gamma}^t-\sum^3_{m=1} \partial_{\xi_m}\Pmp{\xi}D_{x_m}W^t_{\mu,\gamma})\\
&=&-\textrm{Op}(\sum^3_{m=1} \partial_{\xi_m}\Pmp{\xi}D_{x_m}W^t_{\mu,\gamma})
=-\textrm{Op}(\Pmp{D}W^t_{\mu,\gamma}).\\
\end{eqnarray*}
The second equality holds because
\[W_{\mu,\gamma}\Pmp{\xi}-\Pmp{\xi}W_{\mu,\gamma}^t=
\frac{1}{2}\Ppm{D\beta,D\alpha}\Pmp{\xi}-\frac{1}{2}\Pmp{\xi}\Ppm{D\alpha,D\beta}=0,\]
which is true by (\ref{commutation_gen}). Consequently
\begin{eqnarray*}
W_{\mu,\gamma}\Pmp{D}-\Pmp{D}W_{\mu,\gamma}^t&=&-\textrm{Op}(\Pmp{D}W^t_{\mu,\gamma})\\
&=&-\frac{1}{2}\fbM{P_-(D)P_-(D\alpha)}{0}{0}{P_+(D)P_+(D\beta)}-\Pmp{D\kappa}.
\end{eqnarray*}
On the other hand
\begin{equation*}
-W_{\mu,\gamma}W_{\mu,\gamma}^t=-\fbM{(\kappa^2+\frac{1}{4}(D\alpha)^2)I_4}{P_+(D\kappa)}{P_-(D\kappa)}{(\kappa^2+\frac{1}{4}(D\beta)^2)I_4},
\end{equation*}
since
\[\frac{1}{2}(D\alpha+D\beta)=D\textrm{log}(\mu\gamma)^{1/2}=\kappa^{-1}D\kappa.\]
This shows the desired form of $\tilde{Q}_{\mu,\gamma}$.

It remains to prove that $\tilde{Q}'_{\mu,\gamma}$ has the expression given. We have
\[(\Pmp{D}-W_{\mu,\gamma}^t)(\Pmp{D}+W_{\mu,\gamma})=-\Delta I_8-W_{\mu,\gamma}^t\Pmp{D} +\Pmp{D}W_{\mu,\gamma}-W_{\mu,\gamma}^tW_{\mu,\gamma}.\]
Note that $W_{\mu,\gamma}^t=W_{\gamma,\mu}$. Hence, changing the role of $\mu$ and $\gamma$ we obtain
\begin{eqnarray*}
-W_{\mu,\gamma}^t\Pmp{D}+\Pmp{D}W_{\mu,\gamma}&=& -W_{\gamma,\mu}\Pmp{D}+\Pmp{D}W_{\gamma,\mu}^t\\
&=&\frac{1}{2}\fbM{P_-(D)P_-(D\beta)}{0}{0}{P_+(D)P_+(D\alpha)}+\Pmp{D\kappa}.
\end{eqnarray*}
and
\[-W_{\mu,\gamma}^tW_{\mu,\gamma}=-W_{\gamma,\mu}W_{\gamma,\mu}^t=
-\fbM{(\kappa^2+\frac{1}{4}(D\beta)^2)I_4}{P_+(D\kappa)}
{P_-(D\kappa)}{(\kappa^2+\frac{1}{4}(D\alpha)^2)I_4}.\]
This shows the desired form of $\tilde{Q}'_{\mu,\gamma}$.
\end{proof}\\

The crucial point in the above lemma is of course the vanishing of the first order terms. Also, it is worth pointing out that we will need both second order operators given in Lemma \ref{lemma:schroedinger} in our work.

\section{Construction of CGO solutions.}
In this section we construct CGO solutions for the rescaled Maxwell system and then we produce solutions for the original Maxwell system. In order to get such solutions we apply the Sylvester-Uhlmann method (given in \cite{S-U}) to the matrix Schr\"odinger equation obtained in Lemma \ref{lemma:schroedinger}. This argument was introduced in \cite{O-S}.\\

For the sake of completeness we recall the Sylvester-Uhlmann estimate whose proof can be found in \cite{S-U}, \cite{S}. Let $q\in L^\infty(\mathbb{R}^3;\mathbb{C})$ with compact support and $-1<\delta<0$. If $\zeta\in\mathbb{C}^3$ such that $\zeta\cdot\zeta=k^2$ with $|\zeta|$ large enough, then for any $f\in L^2_{\delta+1}$ there exists a unique $u\in H^2_\delta$ solving
\[(\Delta +2\zeta\cdot D +q)u=f.\]
Moreover
\[\norm{u}{}{H^s_\delta}\leq C|\zeta|^{s-1}\norm{f}{}{L^2_{\delta+1}},\]
for $0\leq s\leq 2$ with $C$ a constant independent of $\zeta$ and $f$.
The norms are given by
\[\norm{f}{2}{L^2_{\delta}}=\int_{\mathbb{R}^3}(1+|x|^2)^\delta|f(x)|^2\,dx\]
and
\[
\norm{u}{}{H^s_\delta}=\norm{(1+|x|^2)^{\delta/2}u}{}{H^s(\mathbb{R}^3)}.
\]\\

Assume that we can extend  $\mu$ and $\gamma$ to the whole space so that for some positive constants \(\varepsilon _0\) and \(\mu _0\) we have \(\gamma- \varepsilon _0\), \(\mu - \mu _0\in C^4_0(\R ^3)\),
and let $Q_{\mu,\gamma}=\omega^2\varepsilon_0\mu_0I_8+\tilde{Q}_{\mu,\gamma}$ which is thus compactly supported. Denoting $k=\omega(\varepsilon_0\mu_0)^{1/2}$ we are thus led to the equation
\begin{equation}\label{Schroedinger_eq}
\left(-(\Delta+k^2)I_8+Q_{\mu,\gamma}\right)Z=0
\end{equation}
Let $\theta\in\mathbb{R}^3$ any non-zero vector. Consider $\zeta\in\mathbb{C}^3$ such that $\zeta\cdot\zeta=k^2$ and
\begin{equation*}
\zeta=\tau\hz+\frac{1}{2}\theta+\mathcal{O}(\tau^{-1})
:=\tau\hz+\frac{1}{2}\theta(\tau),
\end{equation*}
where $\tau\geq 1$ is a free parameter controlling the size of $\zeta$, and $\hz=\eta_1+i\eta_2$ with $\eta_1,\eta_2\in\mathbb{R}^3$, $\eta_1\cdot\eta_2=0$ and $|\eta_1|=|\eta_2|=1$. Let $Z_0=Z_0(\zeta)$ be a vector which does not depend on $x$ and which is bounded with respect to the parameter $\tau$.

\begin{prop}\label{proposition:CGO}
Assume that $-1<\delta<0$, and consider $\varepsilon>0$ such that $-1<\delta+\varepsilon<0$. There exists a CGO solution of the form
\[Z=e^{i\zeta\cdot x}(Z_0+Z_{-1}+Z_r)\]
for the equation (\ref{Schroedinger_eq}) in $\mathbb{R}^3$, such that
\begin{gather*}
\norm{D^\alpha Z_{-1}}{}{L^2_{\delta}}=\mathcal{O}(\tau^{-1}),\qquad \norm{D^\alpha Z_r}{}{L^2_{\delta}}=\mathcal{O}(\tau^{|\alpha|-(1+\varepsilon)}),
\end{gather*}
for $0\leq|\alpha|\leq 2$. Moreover, if we denote
\[\hat R=\lim_{\tau\to+\infty}\tau Z_{-1} \qquad \hat M=\lim_{\tau\to+\infty}Z_0,\]
then
\begin{equation}\label{equation_term_-1}
2(\hz\cdot D)I_8 \hat R=-Q_{\mu,\gamma}\hat M
\end{equation}
in $\mathbb{R}^3$.
\end{prop}

\begin{proof}
Let $\zeta\in\mathbb{C}^3$ be as above, note that
\begin{gather*}
e^{-i\zeta\cdot x}\left(-(\Delta+k^2)I_8+Q_{\mu,\gamma}\right)e^{i\zeta\cdot x} (Z_0+Z_{-1}+Z_r)=\\
=Q_{\mu,\gamma}Z_0+((-\Delta+2\zeta\cdot D)I_8+Q_{\mu,\gamma})(Z_{-1}+Z_r)
\end{gather*}
in $\mathbb{R}^3$. Let $\chi$ be a cut-off function which is identically $1$ on a neighborhood of $\{x\in\mathbb{R}^3:|x|<\rho\}$ where \(\rho\) is so large that \(\Omega\subset \{ |x|<\rho\}\), and define
\[
(\hz\cdot D)^{-1}f=\frac{1}{(2\pi)^3}\int_{\mathbb{R}^3}\frac{e^{i\xi\cdot x}\hat f}{\hz\cdot\xi}\,d\xi.
\]
For \(\tau >0\) let $\chi_\tau(x)=\chi(x/\tau)$ and
\[Z_{-1}=\frac{1}{2\tau}\chi_\tau(\hz\cdot D)^{-1}(-Q_{\mu,\gamma}Z_0).\]
It is well known that $(\hz\cdot D)^{-1}:L^2_{\delta+1}\to L^2_{\delta}$ is bounded with norm of order $|\hz|^{-1}$ (\cite{S-U}), hence
\[
\norm{Z_{-1}}{}{L^2_{\delta}}\leq\frac{C}{\tau}\norm{Q_{\mu,\gamma}Z_0}{}{L^2_{\delta+1}}=\mathcal{O}(\tau^{-1}).
\]
Moreover, differentiating we get
\[\norm{D^\alpha Z_{-1}}{}{L^2_{\delta}}=\mathcal{O}(\tau^{-1}),\]
for all $1\leq|\alpha|\leq 2$. Here we need the \(C^4\)-regularity of the coefficients. On the other hand, let $Z_r$ be the solution of
\[
-((-\Delta+2\zeta\cdot D)I_8+Q_{\mu,\gamma})Z_r=((-\Delta+\theta(\tau)\cdot D)I_8+Q_{\mu,\gamma})Z_{-1}+(\hz\cdot D\chi_\tau)(\hz\cdot D)^{-1}(-Q_{\mu,\gamma}Z_0),
\]
satisfying 
\[
\norm{D^\alpha Z_r}{}{L^2_\delta}\leq C|\zeta|^{|\alpha|-1}\norm{((-\Delta+\theta(\tau)\cdot D)I_8+Q_{\mu,\gamma})Z_{-1}+(\hz\cdot D\chi_\tau)(\hz\cdot D)^{-1}(-Q_{\mu,\gamma}Z_0)}{}{L^2_{\delta+1}},
\]
for $0\leq|\alpha|\leq 2$. We then have the estimates
\[\norm{Q_{\mu,\gamma}Z_{-1}}{2}{L^2_{\delta+1}}=\mathcal{O}(\tau^{-2}),\]
\[\norm{(\hz\cdot D\chi_\tau)(\hz\cdot D)^{-1}(-Q_{\mu,\gamma}Z_0)}{2}{L^2_{\delta+1}}\leq C\tau^{-2\varepsilon}\norm{(\hz\cdot D)^{-1}(-Q_{\mu,\gamma}Z_0)}{2}{L^2_{\delta+\varepsilon}}\]
and
\begin{eqnarray*}
\norm{(-\Delta+\theta(\tau)\cdot D)I_8Z_{-1}}{2}{L^2_{\delta+1}}&\leq& \frac{C}{\tau^2}\norm{(\hz\cdot D)^{-1}(-Q_{\mu,\gamma}Z_0)}{2}{L^2_\delta}
+\frac{C}{\tau^2}\norm{D(\hz\cdot D)^{-1}(-Q_{\mu,\gamma}Z_0)}{2}{L^2_\delta}\\
& &+C\tau^{-2\varepsilon}\norm{(-\Delta+\theta(\tau)\cdot D)(\hz\cdot D)^{-1}(-Q_{\mu,\gamma}Z_0)}{2}{L^2_{\delta+\varepsilon}}.
\end{eqnarray*}
Let us explain the second inequality in more detail: consider $\varepsilon>0$ such that $-1<\delta+\varepsilon<0$. Then
\begin{gather*}
\norm{(\hz\cdot D\chi_\tau)(\hz\cdot D)^{-1}(-Q_{\mu,\gamma}Z_0)}{2}{L^2_{\delta+1}}=
\int_{\mathbb{R}^3}(1+|x|^2)^{\delta+1}|(\hz\cdot D\chi_\tau)(\hz\cdot D)^{-1}(-Q_{\mu,\gamma}Z_0)|^2\,dx\\
\leq \sup_{x\in\textrm{supp}(\chi_\tau)}\{(1+|x|^2)^{1-\varepsilon}|\hz\cdot D\chi_\tau|^2\}
\int_{\mathbb{R}^3}(1+|x|^2)^{\delta+\varepsilon}|(\hz\cdot D)^{-1}(-Q_{\mu,\gamma}Z_0)|^2\,dx\\
\leq C\tau^{2(1-\varepsilon)}\tau^{-2}\norm{(\hz\cdot D)^{-1}(-Q_{\mu,\gamma}Z_0)}{2}{L^2_{\delta+\varepsilon}}.
\end{gather*}
Therefore
\[\norm{D^\alpha Z_r}{}{L^2_\delta}=\mathcal{O}(\tau^{|\alpha|-(1+\varepsilon)}),\]
for $0\leq|\alpha|\leq 2$. Finally, $\hat R$ satisfies the equation \eqref{equation_term_-1} by the definition of $Z_{-1}$ and the continuity of $(\hz\cdot D)^{-1}$.
\end{proof}\\

Considering $Z$ as in Proposition \ref{proposition:CGO} we have, by Lemma \ref{lemma:schroedinger}, that 
\[Y=(\Pmp{D}-W^t_{\mu,\gamma})Z\]
is a solution of the rescaled system
\[(\Pmp{D}+W_{\mu,\gamma})Y=0.\]
Note that this solution can be written as
\begin{equation}\label{batman}
Y=e^{i\zeta\cdot x}(Y_1+Y_0+Y_r)
\end{equation}
 where
\begin{equation}\label{robin}
\begin{array}{l}
Y_1=\Pmp{\zeta}Z_0,\\
Y_0=\Pmp{\zeta}Z_{-1}-W^t_{\mu,\gamma}Z_0,\\
Y_r=(\Pmp{D}-W^t_{\mu,\gamma})Z_{-1}+(\Pmp{D+\zeta}-W^t_{\mu,\gamma})Z_r,
\end{array}
\end{equation}
and
\begin{equation}\label{decay_Y}
\norm{Y_1}{}{L^2(U)}=\mathcal{O}(\tau),\quad \norm{D^\alpha Y_0}{}{L^2(U)}=\mathcal{O}(1),\quad \norm{D^\alpha Y_r}{}{L^2(U)}=\mathcal{O}(\tau^{|\alpha|-\varepsilon})
\end{equation}
for $0\leq|\alpha|\leq 1$ and $U$ any bounded open subset of $\mathbb{R}^3$. Recall that our aim is to construct CGO solutions for the original Maxwell system. As we mentioned above, $Y$ will produce a solution of the original Maxwell system if $Y^{\Phi}=Y^{\Psi}=0$. In order to get this condition we choose $Z_0$ in a special way.

\begin{lemma}
If
\[((\Pmp{\zeta}-k)Z_0)^\Phi=((\Pmp{\zeta}-k)Z_0)^\Psi=0\]
then we have
\[Y^t=\fcv{0}{(Y^H)^t}{0}{(Y^E)^t}\]
for $\tau$ large enough.
\end{lemma}
\begin{proof} This proof was given in \cite{O-S}. Here it is included for the sake of completeness.\\
As we mentioned above $Y$ is solution of
\[(\Pmp{D}+W_{\mu,\gamma})Y=0,\]
so it is also solution of
\[(\Pmp{D}-W_{\mu,\gamma}^t)(\Pmp{D}+W_{\mu,\gamma})Y=0.\]
From Lemma \ref{lemma:schroedinger} we know that $Y^\Phi$ and $Y^\Psi$ are solutions of
\begin{equation}
\left\{\begin{array}{l}
-(\Delta+k^2)Y^\Phi+q_\beta Y^\Phi=0\\
-(\Delta+k^2)Y^\Psi+q_\alpha Y^\Psi=0,
\end{array}\right.
\end{equation}
where the potentials
\[q_\beta=-\frac{1}{2}\Delta\beta-\kappa^2+k^2-\frac{1}{4}(D\beta)^2,\qquad q_\alpha=-\frac{1}{2}\Delta\alpha-\kappa^2+k^2-\frac{1}{4}(D\alpha)^2\]
are compactly supported.

Adding $\pm kI_8$ one gets
\[Y=e^{i\zeta\cdot x}(\tilde Y_1+\tilde Y_0+Y_r)\]
where $\tilde Y_1=(\Pmp{\zeta}-kI_8)Z_0$ and $\tilde Y_0=\Pmp{\zeta}Z_{-1}+(kI_8-W^t_{\mu,\gamma})Z_0$. This shows that
\begin{equation}
\left\{\begin{array}{l}
(-\Delta+2\zeta\cdot D+q_\beta)(\tilde Y_1+\tilde Y_0+Y_r)^\Phi=0,\\
(-\Delta+2\zeta\cdot D+q_\alpha)(\tilde Y_1+\tilde Y_0+Y_r)^\Psi=0.
\end{array}\right.
\end{equation}
It is known that these equations have unique solutions in $L^2_\delta$ for $|\zeta|$ large enough (see \cite{S-U}). Since neither $\tilde Y_1^\Phi$ nor $\tilde Y_1^\Psi$ belong to $L^2_\delta$ unless $\tilde Y_1^\Phi=\tilde Y_1^\Psi=0$, we choose $Z_0$ such that
\[\tilde Y_1^\Phi=((\Pmp{\zeta}-kI_8)Z_0)^\Phi=0,\qquad \tilde Y_1^\Psi=((\Pmp{\zeta}-kI_8)Z_0)^\Psi=0.\]
In this way we have $(\tilde Y_1+\tilde Y_0+Y_r)^\Phi,(\tilde Y_1+\tilde Y_0+Y_r)^\Psi\in L^2_\delta$. Then,  the unique solvability for $|\zeta|$ large implies that $Y^\Phi=Y^\Psi=0$ for $|\zeta|$ large enough.
\end{proof}\\

An explicit choice of  $Z_0$ verifying the condition of last lemma may be done as follows. Let 
\setlength\arraycolsep{2pt}
\begin{eqnarray}\label{original_M_form}
Z_0&=&\frac{1}{\tau}\frv{\zeta\cdot a}{kb}{\zeta\cdot b}{ka} =\frac{1}{\tau}\left(\begin{array}{c c|c c}
* & 0 & * & \zeta\cdot\\

* & kI_3 & * & 0\\
\hline
* & \zeta\cdot & * & 0\\

* & 0 & * & kI_3\\
\end{array}
 \right)
\frv{0}{b}{0}{a}=\frac{1}{\tau}\left(\begin{array}{c c|c c}
k & 0 & 0 & \zeta\cdot\\

0 & kI_3 & \zeta & 0\\
\hline
0 & \zeta\cdot & k & 0\\

\zeta & 0 & 0 & kI_3\\
\end{array}
 \right)
\frv{0}{b}{0}{a}\nonumber\\
&=&\frac{1}{\tau}\left(kI_8+\frac{1}{2}(\Ppm{\zeta}+\Pmp{\zeta})\right)m,
\end{eqnarray}
where we denoted $m^t=\fcv{0}{b^t}{0}{a^t}$ for any $a,b\in\mathbb{C}^3$. The same choice was done in \cite{O-S}.

Next we perform the reflection argument given by Isakov in \cite{I}. As in the statement of Theorem 1, let \(\Omega\subset \{x_3 <0\}\) be a $C^{1,1}$ domain. With respect to the standard Cartesian basis $e=\{e_1,e_2,e_3\}$, we introduce the reflection
\[
x=(x_1,x_2,x_3)\mapsto (x_1,x_2,-x_3) =: \dot x (x)
\]
and the reflected domain $\dot \Omega=\{\dot x(x)\in\mathbb{R}^3:x\in\Omega\}$. Consider \[O=\Omega\cup\textrm{int}(\Gamma_0)\cup\dot\Omega.\]
By the assumption \eqref{BoundaryExtensionAssumption}, we can extend the coefficients $\mu_j$ and $\gamma_j$ into $\R^3$ as $C^4$ functions which are even with respect to $x_3$ for $j=1,2$. We still denote the extended coefficients by $\mu_j$ and $\gamma_j$. Note that these functions have the proper regularity to construct the CGO solutions as above.

For the pairs $(\mu_1,\gamma_1)$ and $(\mu_2,\overline{\gamma}_2)$, and for a given vector $\xi = (\xi', \xi_3) \in\mathbb{R}^3$ satisfying $|\xi'| > 0$, we construct CGO solutions $Z_1,Y_1,X_1$ and $Z_2,Y_2,X_2$ in $\R^3$ with complex vectors $\zeta_1$ and $\zeta_2$, respectively. We choose 
\begin{gather}
\zeta_1=\frac{1}{2}\xi+i\left(\tau^2+\frac{|\xi|^2}{4}\right)^{1/2}\eta_1+ \left(\tau^2+k^2\right)^{1/2}\eta_2 \label{zeta1def} \\
\zeta_2=-\frac{1}{2}\xi-i\left(\tau^2+\frac{|\xi|^2}{4}\right)^{1/2}\eta_1+ \left(\tau^2+k^2\right)^{1/2}\eta_2, \label{zeta2def}
\end{gather}
with $\tau\geq 1$ a free parameter controlling the size of $|\zeta_1|$ and $|\zeta_2|$, where $\eta_1,\eta_2\in\mathbb{R}^3$ verify $|\eta_1|=|\eta_2|=1$, $\eta_1\cdot\eta_2=0$ and $\eta_j\cdot\xi=0$ for $j=1,2$. More precisely, for $\xi=(\xi_1,\xi_2,\xi_3)=(\xi',\xi_3)$ we choose
\[\eta_1=\frac{1}{|\xi'|}(\xi_2,-\xi_1,0),\quad \eta_2=\eta_1\wedge\frac{\xi}{|\xi|}=\frac{1}{|\xi'||\xi|} (-\xi_1\xi_3,-\xi_2\xi_3,|\xi'|^2).\]
Note $\zeta_j\cdot\zeta_j=k^2$ for $j=1,2$. Let us denote
\begin{gather*}
\hz=\lim_{\tau\to+\infty}\frac{\zeta_1}{\tau}= \lim_{\tau\to+\infty}\frac{\overline{\zeta}_2}{\tau}=\eta_2+i\eta_1\\
\chz=\lim_{\tau\to+\infty}\frac{\zeta_2}{\tau}=\eta_2-i\eta_1
\end{gather*}
and
\[\xi_1(\tau)=2(\zeta_1-\tau\hz) \qquad \xi_2(\tau)=2(\zeta_2-\tau\chz).\]
Observe that $\xi_1(\tau)=\xi+\mathcal{O}(\tau^{-1})$ and $\xi_2(\tau)=-\xi+\mathcal{O}(\tau^{-1})$. Furthermore, if we denote
\begin{gather*}
\hat R_1=\lim_{\tau\to+\infty}\tau Z^1_{-1} \qquad
\check R_2=\lim_{\tau\to+\infty}\tau Z^2_{-1}\\
\hat M_1=\lim_{\tau\to+\infty}Z^1_0 \qquad \check M_2=\lim_{\tau\to+\infty}Z^2_0
\end{gather*}
then the equations
\begin{gather}
\label{equation_R_1}
2(\hz\cdot D)I_8\hat R_1=-Q_{\mu_1,\gamma_1}\hat M_1\\
\label{equation_R_2}
2(\chz\cdot D)I_8\check R_2=-Q_{\mu_2,\overline{\gamma}_2}\check M_2
\end{gather}
hold in $\mathbb{R}^3$.\\

We now fix another orthonormal basis $f=\{f_1,f_2,f_3\}$ with
\begin{equation} \label{fbasis}
f_2=\frac{1}{|\xi'|}(\xi_1,\xi_2,0),\quad f_3=(0,0,1)=e_3,\quad f_1=f_2\wedge f_3.
\end{equation}
With respect to the basis $f$ one has
$$\xi=(0,|\xi'|,\xi_3)_f,$$
$$\eta_1=(1,0,0)_f=f_1,$$
$$\eta_2=|\xi|^{-1}(0,-\xi_3,|\xi'|)_f.$$
Therefore,
\begin{gather*}
\zeta_1=\left(i\Big(\tau^2+\frac{|\xi|^2}{4}\Big)^{1/2},
\frac{|\xi'|}{2}-(\tau^2+k^2)^{1/2}\frac{\xi_3}{|\xi|},
\frac{\xi_3}{2}+(\tau^2+k^2)^{1/2}\frac{|\xi'|}{|\xi|}\right)_f\\
\zeta_2=\left(-i\Big(\tau^2+\frac{|\xi|^2}{4}\Big)^{1/2},
-\frac{|\xi'|}{2}-(\tau^2+k^2)^{1/2}\frac{\xi_3}{|\xi|},
-\frac{\xi_3}{2}+(\tau^2+k^2)^{1/2}\frac{|\xi'|}{|\xi|}\right)_f.
\end{gather*}
It is obvious that for any $x,y\in\mathbb{R}^3$
\[x\cdot y=\sum^3_{j=1}x_jy_j=\sum^3_{j=1}x^f_jy^f_j\]
where $x=(x_1,x_2,x_3)=(x^f_1,x^f_2,x^f_3)_f$ and $y=(y_1,y_2,y_3)=(y^f_1,y^f_2,y^f_3)_f$ with $x_3=x^f_3$ and $y_3=y^f_3$.

We next show that the reflection argument allows to construct CGO solutions whose tangential magnetic fields $\nu\wedge H$ on $\partial\Omega$ have support contained in $\Gamma$. To this end we introduce
\[\dot X_j(x)=\diag{-\dot I_4}{\dot I_4}X_j(\dot x(x))\]
for $j=1,2$, where
\[\dot I_4=\left(
\begin{array}{c c c c}
1 & 0 & 0 & 0\\
0 & 1 & 0 & 0\\
0 & 0 & 1 & 0\\
0 & 0 & 0 & -1\\
\end{array}
\right).\]
It is an easy matter to check that $\dot X_1$ and $\dot X_2$ are solutions of
\begin{gather*}
(\Pmp{D}+V_{\mu_1,\gamma_1})\dot X_1=0\\
(\Pmp{D}+V_{\mu_2,\overline{\gamma}_2})\dot X_2=0
\end{gather*}
in $\Omega$ (but also in $\dot\Omega$). Therefore, for \(\tau\) large enough,
\[\frv{0}{\mathfrak{H}_1}{0}{\mathfrak{E}_1}=\mathfrak{X}_1=X_1+\dot X_1,\qquad \frv{0}{\mathfrak{H}_2}{0}{\mathfrak{E}_2}=\mathfrak{X}_2=X_2+\dot X_2\]
are solutions of the same equations and $\nu\wedge\mathfrak{H}_j|_{\Gamma_0}=0$ for $j=1,2$. 
Let us state the result obtained as the following proposition.
\begin{prop}\label{reflected_CGO}
Corresponding to pairs $(\mu_1,\gamma_1)$ and $(\mu_2,\overline{\gamma}_2)$, and for a given vector $\xi\in\mathbb{R}^3$ satisfying $|\xi'|>0$, there are CGO solutions $\mathfrak{X}_1$ and $\mathfrak{X}_2$ of the augmented Maxwell equations in $\Omega$ such that
\begin{gather*}
(\Pmp{D}+V_{\mu_1,\gamma_1})\mathfrak{X}_1=0\\
(\Pmp{D}+V_{\mu_2,\overline{\gamma}_2})\mathfrak{X}_2=0
\end{gather*}
where the terms \(\mathfrak{X}_j\) satisfy 
\[\mathfrak{X}_1 = \diag{\mu _1 ^{-1/2}}{\gamma _1^{-1/2}} (Y_1+\diag{-\dot I_4}{\dot I_4}Y_1(\dot x))\]
\[\mathfrak{X}_2 = \diag{\mu _2 ^{-1/2}}{\overline{\gamma}_2^{-1/2}} (Y_2+\diag{-\dot I_4}{\dot I_4}Y_2(\dot x))\]
and the rescaled fields \(Y_j\) are given by \eqref{batman} and \eqref{robin} for complex vectors $\zeta_j$ as in \eqref{zeta1def} and \eqref{zeta2def}, with $\xi = \zeta _1 - \overline{\zeta} _2$. Moreover, for  \(\tau\) large enough the fields \(\mathfrak{X}_j\) will solve Maxwell equations in \(\Omega\) and the tangential components of their magnetic fields will vanish on \(\Gamma _0\).
\end{prop}

\bigskip

\section{The orthogonality identity.}

We assume that $\omega$ is not a resonant frequency, so the Maxwell equations (\ref{Maxwell_system_in_Omega_Intro}) have a unique solution $(E,H)$ for any prescribed value of $\nu\wedge H$ on $\partial\Omega$ as discussed in the introduction. Our uniqueness proofs are based on the following orthogonality relation involving solutions in $\Omega$ and impedance maps on $\partial \Omega$.

\begin{prop}\label{lemma:integral_formula}
Let $(\mu_1,\gamma_1)$ and $(\mu_2,\gamma_2)$ be two pairs of coefficients such that $\omega$ is not a resonant frequency for either. Let $\Lambda_1$ and $\Lambda_2$ be the corresponding impedance maps. Then, for any $X^t_j=\fcv{0}{H^t_j}{0}{E^t_j}$ satisfying
\begin{gather*}
(\Pmp{D}+V_{\mu_1,\gamma_1})X_1=0\\
(\Pmp{D}+V_{\mu_2,\overline{\gamma}_2})X_2=0
\end{gather*}
in $\Omega$
one has
\[\int_{\Omega} X^*_2(V_{\mu_1,\gamma_1}-V_{\mu_2,\gamma_2})X_1\,dx= i\int_{\partial \Omega}(\Lambda_2-\Lambda_1)(\nu\wedge H_1)\cdot \overline{H}_2\,dS.\]
\end{prop}
\begin{proof}
If
\[(u|v)=\int_{\Omega} v^*u\,dx\quad \textrm{ and }\quad
(u|v)_{\partial \Omega}=\int_{\partial \Omega} v^*u\,dS\]
the integration by parts identity
\[(\Pmp{D}X|X')=-i(\Pmp{\nu}X|X')_{\partial \Omega}+(X|\Pmp{D}X')\]
holds.

Since $\Psi_j=\Phi_j=0$ for $j=1,2$ one has that $(X_1|V_{\mu_2,\gamma_2}^*X_2)=(X_1|V_{\mu_2,\overline{\gamma}_2}X_2)$, hence
\setlength\arraycolsep{2pt}
\begin{eqnarray*}
((V_{\mu_1,\gamma_1}-V_{\mu_2,\gamma_2})X_1|X_2)&=&(V_{\mu_1,\gamma_1}X_1|X_2) -(X_1|V_{\mu_2,\gamma_2}^*X_2)\\
&=&(V_{\mu_1,\gamma_1}X_1|X_2)-(X_1|V_{\mu_2,\overline{\gamma}_2}X_2)\\
&=&-(\Pmp{D}X_1|X_2)+(X_1|\Pmp{D}X_2)\\
&=&i(\Pmp{\nu}X_1|X_2)_{\partial \Omega}\\
&=&i\int_{\partial \Omega}\overline{E}_2\cdot(\nu\wedge H_1)-
\overline{H}_2\cdot(\nu\wedge E_1)\,dS
\end{eqnarray*} 
Let $\tilde{X}^t_2=\fcv{0}{\tilde{H}^t_2}{0}{\tilde{E}^t_2}$ be the solution of the boundary value problem
\begin{equation*}
\left\{
\begin{array}{l}
(\Pmp{D}+V_{\mu_2,\gamma_2})\tilde{X}_2=0 \textrm{ in }\Omega\\
\nu\wedge \tilde{H}_2=\nu\wedge H_1 \textrm{ on }\partial \Omega.
\end{array}\right.
\end{equation*}
Since $\tilde X_2$ on the boundary verifies
\[\overline{E}_2\cdot(\nu\wedge H_1)=\overline{E}_2\cdot(\nu\wedge \tilde H_2)\]
we have
\[((V_{\mu_1,\gamma_1}-V_{\mu_2,\gamma_2})X_1|X_2)=i\int_{\partial \Omega}
\overline{E}_2\cdot(\nu\wedge \tilde H_2)-
\overline{H}_2\cdot(\nu\wedge E_1)\,dS.\]
Finally, note that
\setlength\arraycolsep{2pt}
\begin{eqnarray*}
-i(\Pmp{\nu}\tilde X_2|X_2)_{\partial \Omega}&=&(\Pmp{D}\tilde X_2|X_2)-(\tilde X_2|\Pmp{D}X_2)=-(V_{\mu_2,\gamma_2}\tilde X_2|X_2)+(\tilde X_2|V_{\mu_2,\overline{\gamma}_2}X_2)\\
&=&-(V_{\mu_2,\gamma_2}\tilde X_2|X_2)+(\tilde X_2|V_{\mu_2,\gamma_2}^*X_2)= ((V_{\mu_2,\gamma_2}-V_{\mu_2,\gamma_2})\tilde X_2|X_2)=0.\\
\end{eqnarray*}
Here we used again that $\Psi_2=\Phi_2=\tilde\Psi_2=\tilde\Phi_2=0$. From the last identity one gets
\[i\int_{\partial \Omega}\overline{E}_2\cdot(\nu\wedge \tilde H_2)\,dS=i\int_{\partial \Omega}
\overline{H}_2\cdot(\nu\wedge \tilde E_2)\,dS.\]
Therefore
\begin{eqnarray*}
((V_{\mu_1,\gamma_1}-V_{\mu_2,\gamma_2})X_1|X_2)&=&i\int_{\partial \Omega}
\overline{H}_2\cdot(\nu\wedge \tilde E_2)-
\overline{H}_2\cdot(\nu\wedge E_1)\,dS\\
&=&i\int_{\partial \Omega}\overline{H}_2\cdot(\Lambda_2-\Lambda_1)(\nu\wedge H_1)\,dS,
\end{eqnarray*}
which is the desired integral formula.
\end{proof}\\

We next plug the solutions constructed in Proposition \ref{reflected_CGO} into the integral formula given in Proposition \ref{lemma:integral_formula}. Since $\Lambda_1 a=\Lambda_2 a$ on $\Gamma$ for $a\in TH^{1/2}(\partial\Omega)$ with $\textrm{supp}(a)\subset\Gamma$, we have
\begin{eqnarray*}
\int_\Omega \mathfrak{X}^*_2(V_{\mu_1,\gamma_1}-V_{\mu_2,\gamma_2})\mathfrak{X}_1\,dx&=& i\int_{\partial\Omega}(\Lambda_2-\Lambda_1)(\nu\wedge\mathfrak{H}_1)\cdot \overline{\mathfrak{H}}_2\,dS\\
&=&i\int_{\Gamma_0}(\Lambda_2-\Lambda_1)(\nu\wedge\mathfrak{H}_1)\cdot \overline{\mathfrak{H}}_2\,dS.
\end{eqnarray*}
Note that $\nu\wedge\mathfrak{H}_2|_{\Gamma_0}=0$ means that the tangential components of $\mathfrak{H}_2$ on $\Gamma_0$ vanish. On the other hand, $(\Lambda_2-\Lambda_1)(\nu\wedge\mathfrak{H}_1)$ only has tangential components. Hence
\[(\Lambda_2-\Lambda_1)(\nu\wedge\mathfrak{H}_1)\cdot \overline{\mathfrak{H}}_2=0\]
on $\Gamma_0$. Since the boundary term in the integral formula vanishes we get the relation
\begin{equation*}
\int_\Omega \mathfrak{X}^*_2(V_{\mu_1,\gamma_1}-V_{\mu_2,\gamma_2})\mathfrak{X}_1\,dx=0.
\end{equation*}
Following the same notation we write
\[\dot Y_j(x)=\diag{-\dot I_4}{\dot I_4}Y_j(\dot x(x))\]
and we thus have
\begin{gather*}
\mathfrak{X}_1=X_1+\dot X_1=\diag{\mu^{-1/2}_1}{\gamma^{-1/2}_1}(Y_1+\dot Y_1)\\
\mathfrak{X}_2=X_2+\dot X_2=\diag{\mu^{-1/2}_2}{\overline{\gamma}^{-1/2}_2}(Y_2+\dot Y_2).
\end{gather*}
Since for  \(\tau\) large enough the scalar components vanish, i.e. $\mathfrak{X}^\Phi_j=\mathfrak{X}^\Psi_j=0$ for $j=1,2$, writing $V_{\mu_1,\gamma_1}-V_{\mu_2,\gamma_2}$ as in (\ref{potential_V}) gives
\begin{eqnarray*}
\mathfrak{X}^*_2(V_{\mu_1,\gamma_1}-V_{\mu_2,\gamma_2})\mathfrak{X}_1&=& \mathfrak{X}^*_2\diag{\omega(\mu_1-\mu_2)}{\omega(\gamma_1-\gamma_2)}\mathfrak{X}_1\\
&=&(Y_2+\dot Y_2)^*\diag{\tmu}{\tga}(Y_1+\dot Y_1)
\end{eqnarray*}
where
\[\tmu=\omega\frac{\mu_1-\mu_2}{(\mu_1\mu_2)^{1/2}}, \qquad \tga=\omega\frac{\gamma_1-\gamma_2}{(\gamma_1\gamma_2)^{1/2}}.\]
If now we go back to the integral relation, we get
\begin{eqnarray*}
0&=&\int_\Omega \mathfrak{X}^*_2(V_{\mu_1,\gamma_1}-V_{\mu_2,\gamma_2})\mathfrak{X}_1\,dx=
\int_\Omega (Y_2+\dot Y_2)^*\diag{\tmu}{\tga}(Y_1+\dot Y_1)\,dx\\
&=&\int_\Omega Y^*_2\diag{\tmu}{\tga}Y_1\,dx+\int_\Omega \dot Y^*_2\diag{\tmu}{\tga} \dot Y_1\,dx+\int_\Omega Y^*_2\diag{\tmu}{\tga}\dot Y_1\,dx\\
& &+\int_\Omega \dot Y^*_2\diag{\tmu}{\tga}Y_1\,dx.
\end{eqnarray*}
We gather the first two terms and the last two terms (cross terms) in different expressions. Observe that
\setlength\arraycolsep{1pt}
\begin{eqnarray*}
\int_\Omega \dot Y^*_2\diag{\tmu}{\tga} \dot Y_1\,dx&=&
\int_\Omega (Y_2(\dot x))^*\diag{-\dot I_4}{\dot I_4}\diag{\tmu}{\tga}\diag{-\dot I_4}{\dot I_4}Y_1(\dot x)\,dx\\
&=&\int_\Omega (Y_2(\dot x))^*\diag{\tmu}{\tga}Y_1(\dot x)\,dx\\
&=&\int_{\dot\Omega} Y^*_2\diag{\tmu}{\tga}Y_1\,dx
\end{eqnarray*}
hence
\[\int_\Omega Y^*_2\diag{\tmu}{\tga}Y_1\,dx+\int_\Omega \dot Y^*_2\diag{\tmu}{\tga} \dot Y_1\,dx=\int_O Y^*_2\diag{\tmu}{\tga}Y_1\,dx.\]
With a similar argument the cross terms become
\[\int_\Omega Y^*_2\diag{\tmu}{\tga}\dot Y_1\,dx+
\int_\Omega \dot Y^*_2\diag{\tmu}{\tga}Y_1\,dx=
\int_O \dot Y^*_2\diag{\tmu}{\tga}Y_1\,dx.\]
Therefore,
\begin{equation}\label{before_limiting}
0=\int_\Omega \mathfrak{X}^*_2(V_{\mu_1,\gamma_1}-V_{\mu_2,\gamma_2})\mathfrak{X}_1\,dx =\int_O Y^*_2\diag{\tmu}{\tga}Y_1\,dx+\int_O \dot Y^*_2\diag{\tmu}{\tga}Y_1\,dx.
\end{equation}

Below we will integrate by parts in $O$ many times. Even though $O$ may not be a Lipshitz domain, this procedure is easily justified by integrating by parts separately in the $C^{1,1}$ domains $\Omega$ and $\dot{\Omega}$ and by noting that the contribution from $\{x_3=0\}$ vanishes. Now, the terms on the right hand side of \eqref{before_limiting} may be simplified by writing
\begin{equation}
\label{orthogonal_relation}
\int_O (Y^*_2+\dot Y^*_2)\diag{\tmu}{\tga}Y_1\,dx
=\int_O (Y^*_2+\dot Y^*_2)\diag{\tmu}{\tga}(\Pmp{D}-W^t_{\mu_1,\gamma_1})Z_1\,dx.
\end{equation}
Note now that
\[\Pmp{D}(\diag{\tga}{\tmu}Z_1)=\diag{\tmu}{\tga}\Pmp{D}Z_1+(\Pmp{D}\diag{\tga}{\tmu})Z_1\]
and
\begin{gather*}
\int_O (Y^*_2+\dot Y^*_2) \Pmp{D}(\diag{\tga}{\tmu}Z_1)\,dx=\\
=\int_O (\Pmp{D}(Y_2+\dot Y_2))^*\diag{\tga}{\tmu}Z_1\,dx
-i\int_{\partial O} (Y^*_2+\dot Y^*_2)\Pmp{\nu}\diag{\tga}{\tmu}Z_1\,dS
\end{gather*}
integrating by parts. We know by \eqref{BoundaryGammaAssumption} that $\mu_1=\mu_2$ and $\gamma_1=\gamma_2$ on $\Gamma$, and this is valid on $\partial O$ as well. Since $Y_2+\dot Y_2$ is solution for the rescaled system
\[(\Pmp{D}+W_{\mu_2,\overline{\gamma}_2})(Y_2+\dot Y_2)=0\]
in $O$, we have
\[\int_O (Y^*_2+\dot Y^*_2) \Pmp{D}(\diag{\tga}{\tmu}Z_1)\,dx
=-\int_O (W_{\mu_2,\overline{\gamma}_2}(Y_2+\dot Y_2))^*\diag{\tga}{\tmu}Z_1\,dx.\]
From (\ref{orthogonal_relation}) and the above observations we get
\begin{gather}\label{orthogonal_relation_Y2xZ1}
\int_O (Y^*_2+\dot Y^*_2)\diag{\tmu}{\tga}Y_1\,dx=\\
\nonumber
=-\int_O (Y^*_2+\dot Y^*_2)(W^*_{\mu_2,\overline{\gamma}_2}\diag{\tga}{\tmu}+ (\Pmp{D}\diag{\tga}{\tmu})+\diag{\tmu}{\tga}W^t_{\mu_1,\gamma_1})Z_1\,dx.
\end{gather}
It is a straightforward computation to check
\begin{gather*}
W^*_{\mu_2,\overline{\gamma}_2}\diag{\tga}{\tmu}+ (\Pmp{D}\diag{\tga}{\tmu})+\diag{\tmu}{\tga}W^t_{\mu_1,\gamma_1}=\\
\kappa_2\diag{\tga}{\tmu}+\kappa_1\diag{\tmu}{\tga}+\Ppm{D\hga,D\hmu}+\Pmp{D\tga,D\tmu},
\end{gather*}
where
\[\hmu=\omega\frac{\mu_1+\mu_2}{(\mu_1\mu_2)^{1/2}},
\qquad \hga=\omega\frac{\gamma_1+\gamma_2}{(\gamma_1\gamma_2)^{1/2}}.\]
From now on we denote for short
\[U:=\kappa_2\diag{\tga}{\tmu}+\kappa_1\diag{\tmu}{\tga}+\Ppm{D\hga,D\hmu} +\Pmp{D\tga,D\tmu}.\]
Observe that these computations could be written separately for $Y_2$ and $\dot Y_2$, giving
\[\int_O Y^*_2\diag{\tmu}{\tga}Y_1\,dx=-\int_O Y^*_2UZ_1\,dx,\quad
\int_O \dot Y^*_2\diag{\tmu}{\tga}Y_1\,dx=-\int_O \dot Y^*_2UZ_1\,dx.\]
Recall that solutions of the rescaled Maxwell system and the Schr\"odinger system look respectively like
\[Y_2=e^{i\zeta_2\cdot x}(Y^2_1+Y^2_0+Y^2_r), \quad
Z_1=e^{i\zeta_1\cdot x}(Z^1_0+Z^1_{-1}+Z^1_r).\]
Using the estimates in Proposition \ref{proposition:CGO} and in \ref{decay_Y}, we have
\setlength\arraycolsep{2pt}
\begin{eqnarray*}
\lim_{\tau\to+\infty}\int_O Y^*_2\diag{\tmu}{\tga}Y_1\,dx&=&
-\lim_{\tau\to+\infty}
\int_O e^{i(\zeta_1-\overline{\zeta}_2)\cdot x}(Y^2_1)^*UZ^1_0\,dx\\
& &-\lim_{\tau\to+\infty}
\int_O e^{i(\zeta_1-\overline{\zeta}_2)\cdot x}(Y^2_0)^*UZ^1_0\,dx\\
& &-\lim_{\tau\to+\infty}
\int_O e^{i(\zeta_1-\overline{\zeta}_2)\cdot x}(Y^2_1)^*UZ^1_{-1}\,dx.
\end{eqnarray*}
Noting that $\zeta_1-\overline{\zeta}_2=\xi$ and taking the limit in (\ref{before_limiting}), we obtain
\begin{eqnarray}
0&=&-\lim_{\tau\to+\infty}\label{limit_orthogonal_relation_1}
\int_O e^{i\xi\cdot x}(Y^2_1)^*UZ^1_0\,dx\\
& &-\lim_{\tau\to+\infty}\label{limit_orthogonal_relation_2}
\int_O e^{i\xi\cdot x}(Y^2_0)^*UZ^1_0\,dx\\
& &-\lim_{\tau\to+\infty}\label{limit_orthogonal_relation_3}
\int_O e^{i\xi\cdot x}(Y^2_1)^*UZ^1_{-1}\,dx\\
& &-\lim_{\tau\to+\infty}\label{limit_crossoed_terms}
\int_O \dot Y^*_2UZ_1\,dx.
\end{eqnarray}

\section{Some technical computations.}\label{sec:technical}
In this technical section we compute the limit of the terms given in (\ref{limit_orthogonal_relation_1}), (\ref{limit_orthogonal_relation_2}), (\ref{limit_orthogonal_relation_3}) and (\ref{limit_crossoed_terms}).
\begin{lemma}\label{lemma:computation}
The limit in (\ref{limit_orthogonal_relation_1}) is given by
\begin{gather}
\nonumber
\lim_{\tau\to+\infty}
\int_O e^{i\xi\cdot x}(Y^2_1)^*UZ^1_0\,dx=\\
\nonumber
=k((\hz\cdot\overline{b}_2)(\hz\cdot b_1)+ (\hz\cdot\overline{a}_2)(\hz\cdot a_1))\int_O e^{i\xi\cdot x} (\kappa_1+\kappa_2)(\tmu+\tga)\,dx\\
\nonumber
+\omega(\hz\cdot\overline{b}_2)(\hz\cdot b_1)\int_O e^{i\xi\cdot x}
(-\Delta)\Big(\frac{\mu_1}{\mu_2}\Big)^{1/2}\,dx
+\omega(\hz\cdot\overline{a}_2)(\hz\cdot a_1)\int_O e^{i\xi\cdot x}
(-\Delta)\Big(\frac{\gamma_1}{\gamma_2}\Big)^{1/2}\,dx\\
\nonumber
+2k\omega\int_O e^{i\xi\cdot x}
\left[-(\hz\cdot \overline{b}_2)(\hz\wedge a_1)\cdot D\Big(\frac{\mu_2}{\mu_1}\Big)^{1/2}
+(\hz\cdot \overline{a}_2)(\hz\wedge b_1)\cdot D\Big(\frac{\gamma_2}{\gamma_1}\Big)^{1/2}
\right]\,dx\\
\label{limit_computation_1}
+2k\omega\int_O e^{i\xi\cdot x}
\left[-(\hz\cdot b_1)(\hz\wedge\overline{a}_2)\cdot D\Big(\frac{\mu_1}{\mu_2}\Big)^{1/2} +(\hz\cdot a_1)(\hz\wedge\overline{b}_2)\cdot D\Big(\frac{\gamma_1}{\gamma_2}\Big)^{1/2} \right]\,dx.
\end{gather}
\end{lemma}
\begin{proof}
Remember that
\[Z^1_0=\frac{1}{\tau}\left(kI_8+\frac{1}{2}(\Ppm{\zeta_1}+\Pmp{\zeta_1})\right)m_1\]
with
\[m^t_1=\fcv{0}{b^t_1}{0}{a^t_1},\]
and
\[Y^2_1=\Pmp{\zeta_2}Z^2_0.\]
A simple commutation in the middle term after writing the definion of $U$, gives
\begin{gather}
\label{limit_orthogonal_relation_1.1}
\int_O e^{i\xi\cdot x}(Y^2_1)^*UZ^1_0\,dx=
\frac{k}{\tau}\int_O e^{i\xi\cdot x}(\Pmp{\zeta_2}Z^2_0)^* Um_1\,dx\\
\label{limit_orthogonal_relation_1.2}
+\frac{1}{2\tau}\int_O e^{i\xi\cdot x}(Z^2_0)^*\Pmp{\overline{\zeta}_2} (\Ppm{\zeta_1}+\Pmp{\zeta_1})(\kappa_2\diag{\tmu}{\tga}+\kappa_1\diag{\tga}{\tmu})m_1\,dx\\
\label{limit_orthogonal_relation_1.3}
+\frac{1}{2\tau}\int_O e^{i\xi\cdot x}(Z^2_0)^*\Pmp{\overline{\zeta}_2} (\Ppm{D\hga,D\hmu}+\Pmp{D\tga,D\tmu})(\Ppm{\zeta_1}+\Pmp{\zeta_1})m_1\,dx.
\end{gather}
Before computing the limit of (\ref{limit_orthogonal_relation_1.1}) let us write
\begin{equation}\label{Y21_extended}
Y^2_1=\Pmp{\zeta_2}Z^2_0=\frac{1}{\tau}\frv{k\zeta_2\cdot a_2}
{(\zeta_2\cdot b_2)\zeta_2-k\zeta_2\wedge a_2}{k\zeta_2\cdot b_2}
{(\zeta_2\cdot a_2)\zeta_2+k\zeta_2\wedge b_2}.
\end{equation}
It follows from this expression that
\begin{gather}
\nonumber
\lim_{\tau\to+\infty}
\frac{k}{\tau}\int_O e^{i\xi\cdot x} (\Pmp{\zeta_2}Z^2_0)^* Um_1\,dx=\\
\nonumber
=k\int_O e^{i\xi\cdot x} \fcv{0}{(\hz\cdot \overline{b}_2)\hz^t}
{0}{(\hz\cdot \overline{a}_2)\hz^t}(\kappa_2\diag{\tga}{\tmu}+\kappa_1\diag{\tmu}{\tga})
\frv{0}{b_1}{0}{a_1}\,dx\\
\nonumber
+k\int_O e^{i\xi\cdot x} \fcv{0}{(\hz\cdot \overline{b}_2)\hz^t}
{0}{(\hz\cdot \overline{a}_2)\hz^t}(\Ppm{D\hga,D\hmu} +\Pmp{D\tga,D\tmu})
\frv{0}{b_1}{0}{a_1}\,dx\\
\nonumber
=k\int_O e^{i\xi\cdot x}
\left[(\hz\cdot \overline{b}_2)(\hz\cdot b_1)(\kappa_2\tga+\kappa_1\tmu)
+(\hz\cdot \overline{a}_2)(\hz\cdot a_1)(\kappa_2\tmu+\kappa_1\tga)\right]\,dx\\
\label{limit_computation_1.1}
+2k\omega\int_O e^{i\xi\cdot x}
\left[-(\hz\cdot \overline{b}_2)(\hz\wedge a_1)\cdot D\Big(\frac{\mu_2}{\mu_1}\Big)^{1/2}
+(\hz\cdot \overline{a}_2)(\hz\wedge b_1)\cdot D\Big(\frac{\gamma_2}{\gamma_1}\Big)^{1/2}
\right]\,dx.
\end{gather}
To compute the limit of the term in (\ref{limit_orthogonal_relation_1.2}) observe that one has the relation
\begin{gather}
\nonumber
\frac{1}{2\tau}(Z^2_0)^*\Pmp{\overline{\zeta}_2}(\Ppm{\zeta_1}+\Pmp{\zeta_1})=
\frac{k}{2\tau^2}m^*_2\Pmp{\overline{\zeta}_2}(\Ppm{\zeta_1}+\Pmp{\zeta_1})\\
\label{middle_1}
+\frac{1}{4\tau^2}m^*_2(\Ppm{\overline{\zeta}_2}+\Pmp{\overline{\zeta}_2})
\Pmp{\overline{\zeta}_2}(\Ppm{\zeta_1}+\Pmp{\zeta_1}).
\end{gather}
Using the commutation formulas (\ref{pseudo-commutation1}), (\ref{pseudo-commutation2}) and (\ref{commutation}), and the fact that $\zeta_1-\overline{\zeta}_2=\xi$, we can see that
\begin{gather}
\nonumber
(\Ppm{\overline{\zeta}_2}+\Pmp{\overline{\zeta}_2})
\Pmp{\overline{\zeta}_2}(\Ppm{\zeta_1}+\Pmp{\zeta_1})=\\
\nonumber
=(\zeta_1\cdot \zeta_1)\Pmp{\overline{\zeta}_2} -\Ppm{\xi}\Ppm{\zeta_1}\Pmp{\overline{\zeta}_2}
+(\zeta_1\cdot \zeta_1)\Ppm{\overline{\zeta}_2}
-\Pmp{\xi}\Ppm{\overline{\zeta}_2}\Pmp{\zeta_1}\\
\nonumber
+(\overline{\zeta}_2\cdot \overline{\zeta}_2)(\Ppm{\zeta_1}+\Pmp{\zeta_1})\\
\label{middle_2}
=-\Ppm{\xi}\Ppm{\zeta_1}\Pmp{\overline{\zeta}_2}
-\Pmp{\xi}\Ppm{\overline{\zeta}_2}\Pmp{\zeta_1}
+\mathcal{O}(\tau).
\end{gather}
Having in mind formulas (\ref{middle_1}) and (\ref{middle_2}), we compute the limit of (\ref{limit_orthogonal_relation_1.2}):
\begin{gather}
\nonumber
\lim_{\tau\to+\infty}
\frac{1}{2\tau}\int_O e^{i\xi\cdot x}(Z^2_0)^*\Pmp{\overline{\zeta}_2} (\Ppm{\zeta_1}+\Pmp{\zeta_1})(\kappa_2\diag{\tmu}{\tga}+\kappa_1\diag{\tga}{\tmu})m_1\,dx=\\
\nonumber
=\frac{k}{2}\int_O e^{i\xi\cdot x}m^*_2\Pmp{\hz}(\Ppm{\hz}+\Pmp{\hz}) (\kappa_2\diag{\tmu}{\tga}+\kappa_1\diag{\tga}{\tmu})m_1\,dx\\
\nonumber
-\frac{1}{4}\int_O e^{i\xi\cdot x}m^*_2(\Ppm{\xi}+\Pmp{\xi})\Ppm{\hz}\Pmp{\hz} (\kappa_2\diag{\tmu}{\tga}+\kappa_1\diag{\tga}{\tmu})m_1\,dx\\
\nonumber
=\frac{k}{2}\int_O e^{i\xi\cdot x}m^*_2\Pmp{\hz}\Ppm{\hz} (\kappa_2\diag{\tmu}{\tga}+\kappa_1\diag{\tga}{\tmu})m_1\,dx\\
\label{limit_computation_1.2}
=k\int_O e^{i\xi\cdot x}\left[(\hz\cdot \overline{b}_2)(\hz\cdot b_1) (\kappa_2\tmu+\kappa_1\tga)+(\hz\cdot \overline{a}_2)(\hz\cdot a_1) (\kappa_2\tga+\kappa_1\tmu)\right]\,dx.
\end{gather}
Here we used $\Pmp{\hz}\Pmp{\hz}=\hz\cdot\hz=0$,
\[\Ppm{\hz}\Pmp{\hz}=\Pmp{\hz}\Ppm{\hz}=2\fbM{\subM{0}{0}{0}{\hz\otimes\hz}}{0} {0}{\subM{0}{0}{0}{\hz\otimes\hz}}\]
and
\[(\Ppm{\xi}+\Pmp{\xi})\Ppm{\hz}\Pmp{\hz}=0.\]
Finally we study the limit of (\ref{limit_orthogonal_relation_1.3}). Let us write
\[\Pmp{\overline{\zeta}_2}=\tau\Pmp{\hz}+\frac{1}{2}\Pmp{\overline{\xi_2(\tau)}}\]
where
\[\Pmp{\overline{\xi_2(\tau)}}\to-\Pmp{\xi}\textrm{ as }\tau\to+\infty.\]
We have
\begin{gather*}
\frac{1}{2\tau}\int_O e^{i\xi\cdot x}(Z^2_0)^*\Pmp{\overline{\zeta}_2} (\Ppm{D\hga,D\hmu}+\Pmp{D\tga,D\tmu})(\Ppm{\zeta_1}+\Pmp{\zeta_1})m_1\,dx=\\
=\frac{1}{4\tau}\int_O e^{i\xi\cdot x}(Z^2_0)^*\Pmp{\overline{\xi_2(\tau)}} (\Ppm{D\hga,D\hmu}+\Pmp{D\tga,D\tmu})(\Ppm{\zeta_1}+\Pmp{\zeta_1})m_1\,dx\\
+\frac{1}{2}\int_O e^{i\xi\cdot x}(Z^2_0)^*\Pmp{\hz} (\Ppm{D\hga,D\hmu}+\Pmp{D\tga,D\tmu})(\Ppm{\zeta_1}+\Pmp{\zeta_1})m_1\,dx.
\end{gather*}
We split $Z^2_0$ into pieces depending of order of $\tau$,
\begin{equation}\label{Z20_extended}
Z^2_0=\frac{1}{\tau}\frv{\zeta_2\cdot a_2}{0}{\zeta_2\cdot b_2}{0}+
\frac{1}{\tau}\frv{0}{kb_2}{0}{ka_2}:=\frac{1}{\tau}M^2_1+\frac{1}{\tau}M^2_0,
\end{equation}
and we set
\[\check M_2=\lim_{\tau\to+\infty}Z^2_0=\frv{\chz\cdot a_2}{0}{\chz\cdot b_2}{0}.\]
It follows that
\begin{gather*}
\lim_{\tau\to+\infty}
\frac{1}{2\tau}\int_O e^{i\xi\cdot x}(Z^2_0)^*\Pmp{\overline{\zeta}_2} (\Ppm{D\hga,D\hmu}+\Pmp{D\tga,D\tmu})(\Ppm{\zeta_1}+\Pmp{\zeta_1})m_1\,dx=\\
=-\frac{1}{4}\int_O e^{i\xi\cdot x}\check M^*_2\Pmp{\xi} (\Ppm{D\hga,D\hmu}+\Pmp{D\tga,D\tmu})(\Ppm{\hz}+\Pmp{\hz})m_1\,dx\\
+\frac{1}{2}\int_O e^{i\xi\cdot x}(M^2_0)^*\Pmp{\hz} (\Ppm{D\hga,D\hmu}+\Pmp{D\tga,D\tmu})(\Ppm{\hz}+\Pmp{\hz})m_1\,dx\\
+\lim_{\tau\to+\infty}
\frac{1}{2\tau}\int_O e^{i\xi\cdot x}(M^2_1)^*
\Pmp{\hz}(\Ppm{D\hga,D\hmu}+\Pmp{D\tga,D\tmu})(\Ppm{\zeta_1}+\Pmp{\zeta_1})m_1\,dx.
\end{gather*}
Observe that for any $a\in\mathbb{C}^3$ the formula
\begin{gather*}
\Pmp{a}(\Ppm{D\hga,D\hmu}+\Pmp{D\tga,D\tmu})(\Ppm{\hz}+\Pmp{\hz})m_1=\\
=4\omega
\fbM{\subM{a\cdot D(\gamma_1/\gamma_2)^{1/2}}{0}{-a\wedge D(\gamma_1/\gamma_2)^{1/2}}{0}}{0} {0}{\subM{a\cdot D(\mu_1/\mu_2)^{1/2}}{0}{a\wedge D(\mu_1/\mu_2)^{1/2}}{0}}
\frv{\hz\cdot a_1}{0}{\hz\cdot b_1}{0}
\end{gather*}
holds. Using that $\overline{\zeta}_2 = \tau \widehat{\zeta} - \frac{1}{2} \xi + \mathcal{O}(\tau^{-1})$, the limit is
\begin{gather*}
\lim_{\tau\to+\infty}
\frac{1}{2\tau}\int_O e^{i\xi\cdot x}(Z^2_0)^*\Pmp{\overline{\zeta}_2} (\Ppm{D\hga,D\hmu}+\Pmp{D\tga,D\tmu})(\Ppm{\zeta_1}+\Pmp{\zeta_1})m_1\,dx=\\
=-\omega\int_O e^{i\xi\cdot x}
\left[(\hz\cdot\overline{b}_2)(\hz\cdot b_1)(\xi\cdot D)\Big(\frac{\mu_1}{\mu_2}\Big)^{1/2} +(\hz\cdot\overline{a}_2)(\hz\cdot a_1)(\xi\cdot D)\Big(\frac{\gamma_1}{\gamma_2}\Big)^{1/2} \right]\,dx\\
+2k\omega\int_O e^{i\xi\cdot x}
\left[-(\hz\cdot b_1)(\hz\wedge\overline{a}_2)\cdot D\Big(\frac{\mu_1}{\mu_2}\Big)^{1/2} +(\hz\cdot a_1)(\hz\wedge\overline{b}_2)\cdot D\Big(\frac{\gamma_1}{\gamma_2}\Big)^{1/2} \right]\,dx\\
+\lim_{\tau\to+\infty}\frac{2\omega}{\tau}\int_O e^{i\xi\cdot x}
\left[(\overline{\zeta}_2\cdot\overline{b}_2)(\zeta_1\cdot b_1)
(\hz\cdot D)\Big(\frac{\mu_1}{\mu_2}\Big)^{1/2} +(\overline{\zeta}_2\cdot\overline{a}_2)(\zeta_1\cdot a_1)
(\hz\cdot D)\Big(\frac{\gamma_1}{\gamma_2}\Big)^{1/2} \right]\,dx 
\end{gather*}
In the first term on the right, integration by parts gives 
\begin{gather*}
-\omega\int_O e^{i\xi\cdot x}
\left[(\hz\cdot\overline{b}_2)(\hz\cdot b_1)(\xi\cdot D)\Big(\frac{\mu_1}{\mu_2}\Big)^{1/2} +(\hz\cdot\overline{a}_2)(\hz\cdot a_1)(\xi\cdot D)\Big(\frac{\gamma_1}{\gamma_2}\Big)^{1/2} \right]\,dx\\
=\omega\int_O e^{i\xi\cdot x}
\left[(\hz\cdot\overline{b}_2)(\hz\cdot b_1)(-\Delta)\Big(\frac{\mu_1}{\mu_2}\Big)^{1/2} +(\hz\cdot\overline{a}_2)(\hz\cdot a_1)(-\Delta)\Big(\frac{\gamma_1}{\gamma_2}\Big)^{1/2} \right]\,dx\\
+i\omega\int_{\partial O} e^{i\xi\cdot x}
\left[(\hz\cdot\overline{b}_2)(\hz\cdot b_1)(\nu\cdot D)\Big(\frac{\mu_1}{\mu_2}\Big)^{1/2} +(\hz\cdot\overline{a}_2)(\hz\cdot a_1)(\nu\cdot D)\Big(\frac{\gamma_1}{\gamma_2}\Big)^{1/2} \right]\,dS.
\end{gather*}
Since $\partial_\nu^l\mu_1=\partial_\nu^l\mu_2$ and $\partial_\nu^l\gamma_1=\partial_\nu^l\gamma_2$ for $l=0,1$ on $\Gamma$, the boundary term vanishes. On the other hand, for the last term on the right we integrate by parts to obtain 
\begin{eqnarray*}
\int_O e^{i\xi\cdot x}
(\hz\cdot D)\Big(\frac{\mu_1}{\mu_2}\Big)^{1/2}\,dx
&=&\int_O e^{i\xi\cdot x}
(\hz\cdot D)\left[\Big(\frac{\mu_1}{\mu_2}\Big)^{1/2}-1\right]\,dx\\
&=&-\int_O e^{i\xi\cdot x}
(\hz\cdot \xi)\left[\Big(\frac{\mu_1}{\mu_2}\Big)^{1/2}-1\right]\,dx=0
\end{eqnarray*}
and
\begin{eqnarray*}
\int_O e^{i\xi\cdot x}
(\hz\cdot D)\Big(\frac{\gamma_1}{\gamma_2}\Big)^{1/2}\,dx&=&
\int_O e^{i\xi\cdot x}
(\hz\cdot D)\left[\Big(\frac{\gamma_1}{\gamma_2}\Big)^{1/2}-1\right]\,dx=\\
&=&-\int_O e^{i\xi\cdot x}
(\hz\cdot \xi)\left[\Big(\frac{\gamma_1}{\gamma_2}\Big)^{1/2}-1\right]\,dx=0.
\end{eqnarray*}
We introduced $-1$ in order to get rid of the boundary terms in the integrations by parts. Eventually one has
\begin{gather}
\nonumber
\lim_{\tau\to+\infty}
\frac{1}{2\tau}\int_O e^{i\xi\cdot x}(Z^2_0)^*\Pmp{\overline{\zeta}_2} (\Ppm{D\hga,D\hmu}+\Pmp{D\tga,D\tmu})(\Ppm{\zeta_1}+\Pmp{\zeta_1})m_1\,dx=\\
\nonumber
=\omega\int_O e^{i\xi\cdot x}
\left[(\hz\cdot\overline{b}_2)(\hz\cdot b_1)(-\Delta)\Big(\frac{\mu_1}{\mu_2}\Big)^{1/2} +(\hz\cdot\overline{a}_2)(\hz\cdot a_1)(-\Delta)\Big(\frac{\gamma_1}{\gamma_2}\Big)^{1/2} \right]\,dx\\
\label{limit_computation_1.3}
+2k\omega\int_O e^{i\xi\cdot x}
\left[-(\hz\cdot b_1)(\hz\wedge\overline{a}_2)\cdot D\Big(\frac{\mu_1}{\mu_2}\Big)^{1/2} +(\hz\cdot a_1)(\hz\wedge\overline{b}_2)\cdot D\Big(\frac{\gamma_1}{\gamma_2}\Big)^{1/2} \right]\,dx.
\end{gather}
Putting (\ref{limit_computation_1.1}), (\ref{limit_computation_1.2}) and (\ref{limit_computation_1.3}) together we get the final result.
\end{proof}

\begin{lemma}
The limit in (\ref{limit_orthogonal_relation_2}) is given by
\begin{gather}
\nonumber
\lim_{\tau\to+\infty}\int_O e^{i\xi\cdot x}(Y^2_0)^*UZ^1_0 \,dx=\\
\nonumber
=-\int_O e^{i\xi\cdot x}\kappa_2\left[(\hz\cdot \overline{a}_2)(\hz\cdot a_1) (\kappa_2\tga+\kappa_1\tmu)+(\hz\cdot \overline{b}_2)(\hz\cdot b_1) (\kappa_2\tmu+\kappa_1\tga)\right]\,dx\\
\nonumber
+\omega\int_O e^{i\xi\cdot x}\left[(\hz\cdot \overline{a}_2)(\hz\cdot a_1) \frac{D\gamma_2}{\gamma_2}\cdot D\Big(\frac{\gamma_1}{\gamma_2}\Big)^{1/2}+
(\hz\cdot \overline{b}_2)(\hz\cdot b_1)
\frac{D\mu_2}{\mu_2}\cdot D\Big(\frac{\mu_1}{\mu_2}\Big)^{1/2}\right]\,dx\\
\nonumber
+\int_O e^{i\xi\cdot x}
\left[(\hz\cdot\hat R^E_2)(\hz\cdot a_1) (\kappa_2\tga+\kappa_1\tmu)+
(\hz\cdot\hat R^H_2)(\hz\cdot b_1) (\kappa_2\tmu+\kappa_1\tga)\right]\,dx\\
\nonumber
+2\omega\int_O e^{i\xi\cdot x}\left[(\hz\cdot b_1)(\hat R^E_2\wedge\hz)
\cdot D\Big(\frac{\mu_1}{\mu_2}\Big)^{1/2}
-(\hz\cdot a_1) (\hat R^H_2\wedge\hz)\cdot D\Big(\frac{\gamma_1}{\gamma_2}\Big)^{1/2}\right]\,dx\\
\nonumber
-\omega(\hz\cdot b_1)\int_O e^{i\xi\cdot x}
2(\hz\cdot D)\hat R^\Psi_2\left[\Big(\frac{\mu_1}{\mu_2}\Big)^{1/2}-1\right]\,dx\\
\label{limit_computation_2}
-\omega(\hz\cdot a_1)\int_O e^{i\xi\cdot x}
2(\hz\cdot D)\hat R^\Phi_2 \left[\Big(\frac{\gamma_1}{\gamma_2}\Big)^{1/2}-1\right]\,dx,
\end{gather}
where
\[\hat R^t_2=\fcv{\hat R^\Phi_2}{(\hat R^H_2)^t}{\hat R^\Psi_2}{(\hat R^E_2)^t}=\check R^*_2\]
with
\[\check R_2=\lim_{\tau\to+\infty}\tau Z^2_{-1}.\]
\end{lemma}
\begin{proof}
Since $Y^2_0=\Pmp{\zeta_2}Z^2_{-1}-W^t_{\mu_2,\overline{\gamma}_2}Z^2_0$ we have
\begin{eqnarray}
\label{limit_orthogonal_relation_2.1}
\int_O e^{i\xi\cdot x}(Y^2_0)^*UZ^1_0\,dx&=&
\int_O e^{i\xi\cdot x}
(Z^2_{-1})^*\Pmp{\overline{\zeta}_2}UZ^1_0\,dx\\
\label{limit_orthogonal_relation_2.2}
& &-\int_O e^{i\xi\cdot x}
(Z^2_0)^*\overline{W_{\mu_2,\overline{\gamma}_2}}UZ^1_0\,dx.
\end{eqnarray}
We can write
\[\overline{W_{\mu_2,\overline{\gamma}_2}}=\kappa_2I_8-\frac{1}{2}\Ppm{D\beta_2,D\alpha_2}\]
and set
\[\hat M_1=\lim_{\tau\to+\infty}Z^1_0=\frv{\hz\cdot a_1}{0}{\hz\cdot b_1}{0}
\quad \check M_2=\lim_{\tau\to+\infty}Z^2_0=\frv{\chz\cdot a_2}{0}{\chz\cdot b_2}{0}.\]
Then the limit of (\ref{limit_orthogonal_relation_2.2}) is
\begin{gather}
\nonumber
\lim_{\tau\to+\infty}-\int_O e^{i\xi\cdot x}
(Z^2_0)^*\overline{W_{\mu_2,\overline{\gamma}_2}}UZ^1_0\,dx=\\
\nonumber
=-\int_O e^{i\xi\cdot x}\kappa_2\check M^*_2U\hat M_1\,dx
+\frac{1}{2}\int_O e^{i\xi\cdot x}\check M^*_2\Ppm{D\beta_2,D\alpha_2}U\hat M_1\,dx\\
\nonumber
=-\int_O e^{i\xi\cdot x}\kappa_2\check M^*_2 (\kappa_2\diag{\tga}{\tmu}+\kappa_1\diag{\tmu}{\tga})\hat M_1\,dx\\
\nonumber
+\frac{1}{2}\int_O e^{i\xi\cdot x}\check M^*_2\Ppm{D\beta_2,D\alpha_2}
(\Ppm{D\hga,D\hmu} +\Pmp{D\tga,D\tmu})\hat M_1\,dx\\
\nonumber
=-\int_O e^{i\xi\cdot x}\kappa_2\left[(\hz\cdot \overline{a}_2)(\hz\cdot a_1) (\kappa_2\tga+\kappa_1\tmu)+(\hz\cdot \overline{b}_2)(\hz\cdot b_1) (\kappa_2\tmu+\kappa_1\tga)\right]\,dx\\
\label{limit_computation_2.2}
+\omega\int_O e^{i\xi\cdot x}\left[(\hz\cdot \overline{a}_2)(\hz\cdot a_1) \frac{D\gamma_2}{\gamma_2}\cdot D\Big(\frac{\gamma_1}{\gamma_2}\Big)^{1/2}+
(\hz\cdot \overline{b}_2)(\hz\cdot b_1)
\frac{D\mu_2}{\mu_2}\cdot D\Big(\frac{\mu_1}{\mu_2}\Big)^{1/2}\right]\,dx.
\end{gather}
In the second equality we just wrote down $U$ and we used
\[\check M^*_2(\Ppm{D\hga,D\hmu}+\Pmp{D\tga,D\tmu})\hat M_1=0\]
and
\[\check M^*_2\Ppm{D\beta_2,D\alpha_2}(\kappa_2\diag{\tga}{\tmu}+\kappa_1\diag{\tmu}{\tga}) \hat M_1=0.\]
We now compute the limit of the right hand side of (\ref{limit_orthogonal_relation_2.1}):
\begin{gather*}
\lim_{\tau\to+\infty}\int_O e^{i\xi\cdot x}
(Z^2_{-1})^*\Pmp{\overline{\zeta}_2}UZ^1_0\,dx=\\
=\int_O e^{i\xi\cdot x}
\left[(\hz\cdot\hat R^E_2)(\hz\cdot a_1) (\kappa_2\tga+\kappa_1\tmu)+
(\hz\cdot\hat R^H_2)(\hz\cdot b_1) (\kappa_2\tmu+\kappa_1\tga)\right]\,dx\\
+2\omega\int_O e^{i\xi\cdot x}\left[(\hz\cdot b_1)(\hat R^E_2\wedge\hz)
\cdot D\Big(\frac{\mu_1}{\mu_2}\Big)^{1/2}
-(\hz\cdot a_1) (\hat R^H_2\wedge\hz)\cdot D\Big(\frac{\gamma_1}{\gamma_2}\Big)^{1/2}\right]\,dx\\
+2\omega\int_O e^{i\xi\cdot x}\left[(\hz\cdot b_1)
\hat R^\Psi_2(\hz\cdot D)\Big(\frac{\mu_1}{\mu_2}\Big)^{1/2}
+(\hz\cdot a_1) \hat R^\Phi_2(\hz\cdot D)
\Big(\frac{\gamma_1}{\gamma_2}\Big)^{1/2}\right]\,dx.
\end{gather*}
Integrating by parts one gets
\begin{eqnarray*}
\int_O e^{i\xi\cdot x}
\hat R^\Psi_2(\hz\cdot D)\Big(\frac{\mu_1}{\mu_2}\Big)^{1/2}\,dx&=&
\int_O e^{i\xi\cdot x}
\hat R^\Psi_2(\hz\cdot D)\left[\Big(\frac{\mu_1}{\mu_2}\Big)^{1/2}-1\right]\,dx\\
&=&-\int_O e^{i\xi\cdot x}
(\hz\cdot D)\hat R^\Psi_2\left[\Big(\frac{\mu_1}{\mu_2}\Big)^{1/2}-1\right]\,dx
\end{eqnarray*}
and
\begin{eqnarray*}
\int_O e^{i\xi\cdot x}
\hat R^\Phi_2 (\hz\cdot D)\Big(\frac{\gamma_1}{\gamma_2}\Big)^{1/2}\,dx&=&
\int_O e^{i\xi\cdot x}
\hat R^\Phi_2 (\hz\cdot D)\left[\Big(\frac{\gamma_1}{\gamma_2}\Big)^{1/2}-1\right]\,dx\\
&=&-\int_O e^{i\xi\cdot x}
(\hz\cdot D)\hat R^\Phi_2 \left[\Big(\frac{\gamma_1}{\gamma_2}\Big)^{1/2}-1\right]\,dx.
\end{eqnarray*}
Once again, we introduced $-1$ in order to get rid of the boundary terms. Finally
\begin{gather}
\nonumber
\lim_{\tau\to+\infty}\int_O e^{i\xi\cdot x}
(Z^2_{-1})^*\Pmp{\overline{\zeta}_2}UZ^1_0\,dx=\\
\nonumber
=\int_O e^{i\xi\cdot x}
\left[(\hz\cdot\hat R^E_2)(\hz\cdot a_1) (\kappa_2\tga+\kappa_1\tmu)+
(\hz\cdot\hat R^H_2)(\hz\cdot b_1) (\kappa_2\tmu+\kappa_1\tga)\right]\,dx\\
\nonumber
+2\omega\int_O e^{i\xi\cdot x}\left[(\hz\cdot b_1)(\hat R^E_2\wedge\hz)
\cdot D\Big(\frac{\mu_1}{\mu_2}\Big)^{1/2}
-(\hz\cdot a_1) (\hat R^H_2\wedge\hz)\cdot D\Big(\frac{\gamma_1}{\gamma_2}\Big)^{1/2}\right]\,dx\\
\nonumber
-2\omega(\hz\cdot b_1)\int_O e^{i\xi\cdot x}
(\hz\cdot D)\hat R^\Psi_2\left[\Big(\frac{\mu_1}{\mu_2}\Big)^{1/2}-1\right]\,dx\\
\label{limit_computation_2.1}
-2\omega(\hz\cdot a_1)\int_O e^{i\xi\cdot x}
(\hz\cdot D)\hat R^\Phi_2 \left[\Big(\frac{\gamma_1}{\gamma_2}\Big)^{1/2}-1\right]\,dx.
\end{gather}
Putting (\ref{limit_computation_2.2}) and (\ref{limit_computation_2.1}) together we get the final formula.
\end{proof}

\begin{lemma}
The limit in (\ref{limit_orthogonal_relation_3}) is given by
\begin{gather}
\nonumber
\lim_{\tau\to+\infty}
\int_O e^{i\xi\cdot x}(Y^2_1)^*UZ^1_{-1}\,dx=\\
\nonumber
=\int_O e^{i\xi\cdot x}
\left[(\kappa_2\tga+\kappa_1\tmu)(\hz\cdot\overline{b}_2)(\hz\cdot\hat R^H_1)
+(\kappa_2\tmu+\kappa_1\tga)(\hz\cdot\overline{a}_2)(\hz\cdot\hat R^E_1)\right]\,dx\\
\nonumber
+2\omega\int_O e^{i\xi\cdot x}\left[
(\hz\cdot\overline{a}_2)(\hz\wedge\hat R^H_1)\cdot D\Big(\frac{\gamma_2}{\gamma_1}\Big)^{1/2}
-(\hz\cdot\overline{b}_2)(\hz\wedge\hat R^E_1)\cdot D\Big(\frac{\mu_2}{\mu_1}\Big)^{1/2} \right]\,dx\\
\nonumber
-\omega(\hz\cdot\overline{b}_2)\int_O e^{i\xi\cdot x}
2(\hz\cdot D)\hat R^\Psi_1\left[\Big(\frac{\mu_1}{\mu_2}\Big)^{1/2}-1\right]\,dx\\
\label{limit_computation_3}
-\omega(\hz\cdot\overline{a}_2)\int_O e^{i\xi\cdot x}
2(\hz\cdot D)\hat R^\Phi_1 \left[\Big(\frac{\gamma_1}{\gamma_2}\Big)^{1/2}-1\right]\,dx,
\end{gather}
where
\[\hat R^t_1=\fcv{\hat R^\Phi_1}{(\hat R^H_1)^t}{\hat R^\Psi_1}{(\hat R^E_1)^t}\]
with
\[\hat R_1=\lim_{\tau\to+\infty}\tau Z^1_{-1}.\]
\end{lemma}
\begin{proof}
We use the computations done in (\ref{Y21_extended}) for $Y^2_1 = \Pmp{\zeta_2} Z^2_0$ and look at the limit
\begin{gather}
\nonumber
\lim_{\tau\to+\infty}
\int_O e^{i\xi\cdot x}(Y^2_1)^*UZ^1_{-1}\,dx 
=\lim_{\tau\to+\infty}
\int_O e^{i\xi\cdot x} (\Pmp{\zeta_2} Z^2_0)^* UZ^1_{-1}\,dx\\
\label{limit_orthogonal_relation_3.1}
=\int_O e^{i\xi\cdot x}
\fcv{0}{(\hz\cdot \overline{b}_2)\hz^t}
{0}{(\hz\cdot \overline{a}_2)\hz^t}U\hat R_1\,dx.
\end{gather}
We insert the definition of $U$ in (\ref{limit_orthogonal_relation_3.1}) and we compute
\begin{gather*}
\int_O e^{i\xi\cdot x}
\fcv{0}{(\hz\cdot \overline{b}_2)\hz^t}
{0}{(\hz\cdot \overline{a}_2)\hz^t}U\hat R_1\,dx=\\
=\int_O e^{i\xi\cdot x}
\left[(\kappa_2\tga+\kappa_1\tmu)(\hz\cdot\overline{b}_2)(\hz\cdot\hat R^H_1)
+(\kappa_2\tmu+\kappa_1\tga)(\hz\cdot\overline{a}_2)(\hz\cdot\hat R^E_1)\right]\,dx\\
+\int_O e^{i\xi\cdot x}
\fcv{0}{(\hz\cdot \overline{b}_2)\hz^t}
{0}{(\hz\cdot \overline{a}_2)\hz^t}(\Ppm{D\hga,D\hmu}+\Pmp{D\tga,D\tmu})\hat R_1\,dx\\
=\int_O e^{i\xi\cdot x}
\left[(\kappa_2\tga+\kappa_1\tmu)(\hz\cdot\overline{b}_2)(\hz\cdot\hat R^H_1)
+(\kappa_2\tmu+\kappa_1\tga)(\hz\cdot\overline{a}_2)(\hz\cdot\hat R^E_1)\right]\,dx\\
+2\omega\int_O e^{i\xi\cdot x}
\left[\hat R^\Psi_1(\hz\cdot\overline{b}_2)(\hz\cdot D)\Big(\frac{\mu_1}{\mu_2}\Big)^{1/2}
-(\hz\cdot\overline{b}_2)(\hz\wedge\hat R^E_1)\cdot D\Big(\frac{\mu_2}{\mu_1}\Big)^{1/2} \right]\,dx\\
+2\omega\int_O e^{i\xi\cdot x}\left[
\hat R^\Phi_1(\hz\cdot\overline{a}_2)(\hz\cdot D)\Big(\frac{\gamma_1}{\gamma_2}\Big)^{1/2}+ (\hz\cdot\overline{a}_2)(\hz\wedge\hat R^H_1)\cdot D\Big(\frac{\gamma_2}{\gamma_1}\Big)^{1/2} \right]\,dx.
\end{gather*}
Integrating by parts one gets
\begin{eqnarray*}
\int_O e^{i\xi\cdot x}
\hat R^\Psi_1(\hz\cdot D)\Big(\frac{\mu_1}{\mu_2}\Big)^{1/2}\,dx&=&
\int_O e^{i\xi\cdot x}
\hat R^\Psi_1(\hz\cdot D)\left[\Big(\frac{\mu_1}{\mu_2}\Big)^{1/2}-1\right]\,dx\\
&=&-\int_O e^{i\xi\cdot x}
(\hz\cdot D)\hat R^\Psi_1\left[\Big(\frac{\mu_1}{\mu_2}\Big)^{1/2}-1\right]\,dx
\end{eqnarray*}
and
\begin{eqnarray*}
\int_O e^{i\xi\cdot x}
\hat R^\Phi_1 (\hz\cdot D)\Big(\frac{\gamma_1}{\gamma_2}\Big)^{1/2}\,dx&=&
\int_O e^{i\xi\cdot x}
\hat R^\Phi_1 (\hz\cdot D)\left[\Big(\frac{\gamma_1}{\gamma_2}\Big)^{1/2}-1\right]\,dx\\
&=&-\int_O e^{i\xi\cdot x}
(\hz\cdot D)\hat R^\Phi_1 \left[\Big(\frac{\gamma_1}{\gamma_2}\Big)^{1/2}-1\right]\,dx.
\end{eqnarray*}
This gives the final result.
\end{proof}

\medskip
The last technical lemma shows that the cross terms do not contribute to the limit.

\begin{lemma} \label{lemma:computation4}
The limit in (\ref{limit_crossoed_terms}) satisfies 
\begin{equation}\label{limit_computation_crossoed_terms}
\lim_{\tau\to+\infty}\int_O \dot Y^*_2UZ_1\,dx=0.
\end{equation}
\end{lemma}
\begin{proof}
Recall that
\[\dot Y^*_2=e^{-i\overline{\zeta}_2\cdot \dot x}(\dot Y^2_1+\dot Y^2_0+\dot Y^2_r)^*\qquad
Z_1=e^{i\zeta_1\cdot x}(Z^1_0+Z^1_{-1}+Z^1_r),\]
where $\dot Y^2_j=\diag{-\dot I_4}{\dot I_4}Y^2_j(\dot x)$ for $j=1,0,r$. Note that writing the dot product in terms of the basis $f$ in \eqref{fbasis} one has
\[\phi_\tau(x):=\zeta_1\cdot x -\overline{\zeta}_2\cdot \dot x(x)= |\xi'|x^f_2+2(\tau^2+k^2)^{1/2}\frac{|\xi'|}{|\xi|}x^f_3.\]
Since we chose $f_3=e_3$ we have that $x_3^f=x_3$, hence
\begin{gather*}
\int_O \dot Y^*_2UZ_1\,dx=\int_O e^{i\phi_\tau}(\dot Y^2_1+\dot Y^2_0+\dot Y^2_r)^*U(Z^1_0+Z^1_{-1}+Z^1_r)\,dx\\
=\Big(2(\tau^2+k^2)^{1/2}\frac{|\xi'|}{|\xi|}\Big)^{-1} \int_O D_{x_3}e^{i\phi_\tau}(\dot Y^2_1+\dot Y^2_0+\dot Y^2_r)^*U(Z^1_0+Z^1_{-1}+Z^1_r)\,dx\\
=-\Big(2(\tau^2+k^2)^{1/2}\frac{|\xi'|}{|\xi|}\Big)^{-1} \int_O e^{i\phi_\tau}D_{x_3}[(\dot Y^2_1+\dot Y^2_0+\dot Y^2_r)^*U(Z^1_0+Z^1_{-1}+Z^1_r)]\,dx\\
-i\Big(2(\tau^2+k^2)^{1/2}\frac{|\xi'|}{|\xi|}\Big)^{-1} \int_{\partial O} e^{i\phi_\tau}\nu_3(\dot Y^2_1+\dot Y^2_0+\dot Y^2_r)^*U(Z^1_0+Z^1_{-1}+Z^1_r)\,dx.\\
\end{gather*}
The boundary term vanishes because we assumed that $\partial_\nu^l\mu_1=\partial_\nu^l\mu_2$ and $\partial_\nu^l\gamma_1=\partial_\nu^l\gamma_2$ for $l=0,1$ on $\Gamma$. Since $\dot Y^2_1$ and $Z^1_0$ are constant vectors with respect to $x$ we have
\begin{eqnarray*}
D_{x_3}[(\dot Y^2_1+\dot Y^2_0+\dot Y^2_r)^*U(Z^1_0+Z^1_{-1}+Z^1_r)]&=&
D_{x_3}(\dot Y^2_0+\dot Y^2_r)^*U(Z^1_0+Z^1_{-1}+Z^1_r)\\
& &+(\dot Y^2_1+\dot Y^2_0+\dot Y^2_r)^*D_{x_3}U(Z^1_0+Z^1_{-1}+Z^1_r)\\
& &+(\dot Y^2_1+\dot Y^2_0+\dot Y^2_r)^*UD_{x_3}(Z^1_{-1}+Z^1_r).
\end{eqnarray*}
Note that from the asymptotic behavior in (\ref{decay_Y}) and in Proposition \ref{proposition:CGO} one can see that
\begin{gather*}
D_{x_3}(\dot Y^2_0+\dot Y^2_r)^*U(Z^1_0+Z^1_{-1}+Z^1_r)=D_{x_3}(\dot Y^2_r)^*UZ^1_0+\mathcal{O}(1),\\
(\dot Y^2_1+\dot Y^2_0+\dot Y^2_r)^*D_{x_3}U(Z^1_0+Z^1_{-1}+Z^1_r)=(\dot Y^2_1)^*D_{x_3}UZ^1_0+\mathcal{O}(1),\\
(\dot Y^2_1+\dot Y^2_0+\dot Y^2_r)^*UD_{x_3}(Z^1_{-1}+Z^1_r)=(\dot Y^2_1)^*UD_{x_3}Z^1_r+\mathcal{O}(1),
\end{gather*}
in $L^2(O)$; moreover, \[\norm{D_{x_3}(\dot Y^2_r)^*UZ^1_0}{}{L^2(O)}=\mathcal{O}(\tau^{1-\varepsilon})\qquad \norm{(\dot Y^2_1)^*UD_{x_3}Z^1_r}{}{L^2(O)}=\mathcal{O}(\tau^{1-\varepsilon}).\]
On the other hand, we know by (\ref{Y21_extended}) and (\ref{Z20_extended}) that
\begin{gather*}
(Y^2_1)^t=\frac{1}{\tau}\fcv{0}{(\zeta_2\cdot b_2)\zeta^t_2}{0}{(\zeta_2\cdot a_2)\zeta^t_2}+\mathcal{O}(1)
:=\frac{1}{\tau}(L^2_2)^t+\mathcal{O}(1),\\
(Z^1_0)^t=\frac{1}{\tau}\fcv{\zeta_1\cdot a_1}{0}{\zeta_1\cdot b_1}{0}+\mathcal{O}(\tau^{-1})
:=\frac{1}{\tau}(M^1_1)^t+\mathcal{O}(\tau^{-1})
\end{gather*}
in $L^2(O)$. Therefore
\begin{eqnarray*}
\lim_{\tau\to+\infty}\int_O \dot Y^*_2UZ_1\,dx&=&-\lim_{\tau\to+\infty} \Big(2(\tau^2+k^2)^{1/2}\frac{|\xi'|}{|\xi|}\Big)^{-1}
\left[ \int_O e^{i\phi_\tau}D_{x_3}(\dot Y^2_r)^*UZ^1_0\,dx\right.\\
& &\left. +\int_O e^{i\phi_\tau}(\dot Y^2_1)^*D_{x_3}UZ^1_0\,dx+
\int_O e^{i\phi_\tau}(\dot Y^2_1)^*UD_{x_3}Z^1_r\,dx\right]
\end{eqnarray*}
\begin{gather*}
=-\lim_{\tau\to+\infty}\Big(2(\tau^2+k^2)^{1/2}\frac{|\xi'|}{|\xi|}\Big)^{-1}
\int_O e^{i\phi_\tau}(Y^2_1(\dot{x}))^*\diag{-\dot I_4}{\dot I_4}D_{x_3}UZ^1_0\,dx\\
=-\lim_{\tau\to+\infty}\frac{1}{\tau^2}\Big(2(\tau^2+k^2)^{1/2}\frac{|\xi'|}{|\xi|}\Big)^{-1}
\int_O e^{i\phi_\tau}(L^2_2)^*\diag{-\dot I_4}{\dot I_4}D_{x_3}UM^1_1\,dx.
\end{gather*}
Note that
\begin{gather*}
(L^2_2)^*\diag{-\dot I_4}{\dot I_4}D_{x_3}UM^1_1=\\
=2\omega(\overline{\zeta}_2\cdot \overline{a}_2)(\zeta_1\cdot a_1)(\dot{\overline{\zeta}}_2\cdot D) D_{x_3}\Big(\frac{\gamma_1}{\gamma_2}\Big)^{1/2}
-2\omega(\overline{\zeta}_2\cdot \overline{b}_2)(\zeta_1\cdot b_1)(\dot{\overline{\zeta}}_2\cdot D) D_{x_3}\Big(\frac{\mu_1}{\mu_2}\Big)^{1/2}.
\end{gather*}
Since the Riemann-Lebesgue lemma implies that
\[\lim_{\tau\to+\infty}\int_O e^{i\phi_\tau} D_jD_{x_3}\Big(\frac{\mu_1}{\mu_2}\Big)^{1/2}\,dx=
\lim_{\tau\to+\infty}\int_O e^{i\phi_\tau} D_jD_{x_3}\Big(\frac{\gamma_1}{\gamma_2}\Big)^{1/2}\,dx=0,\]
we have proved the result.
\end{proof}

\section{Proof of Theorem 1.}
Assuming the conditions in Theorem 1, we have proved that the identity in  \eqref{limit_orthogonal_relation_1}--\eqref{limit_crossoed_terms} holds and that the terms in that identity have the limits given in Lemmas \ref{lemma:computation} to \ref{lemma:computation4}. In this section we show the equations verified by the coefficients and then give a proof of Theorem 1.
\begin{prop}
Let $\xi \in \R^3$. With appropriate choices of the constant vectors $a_1$, $a_2$, $b_1$ and $b_2$ one has the following equations:
\begin{itemize}
\item [(a)] If we choose $b_1=\overline{b}_2=\chz$ and $a_1=\overline{a}_2=\hz$
\begin{equation}\label{magnetic_laplacian}
\int_O e^{i\xi\cdot x} \left[\frac{1}{2}\Delta(\beta_2-\beta_1)+ \frac{1}{4}[(D\beta_1)^2-(D\beta_2)^2]+\kappa^2_1-\kappa^2_2\right]\,dx=0.
\end{equation}
\item[(b)] If we choose $a_1=\overline{a}_2=\chz$ and $b_1=\overline{b}_2=\hz$
\begin{eqnarray}\label{electric_laplacian}
\int_O e^{i\xi\cdot x} \left[\frac{1}{2}\Delta(\alpha_2-\alpha_1)+ \frac{1}{4}[(D\alpha_1)^2-(D\alpha_2)^2]+\kappa^2_1-\kappa^2_2\right]\,dx=0.
\end{eqnarray}
\end{itemize}
\end{prop}
\begin{proof}
In order to get these equations from the previous lemmas we will need the equations (\ref{equation_R_1}) and (\ref{equation_R_2}). Let us write these equations by components:
\begin{equation}\label{equation_R_1_components}
\left\{\begin{array}{l}
2(\hz\cdot D)\hat R^\Phi_1= -\Big(\frac{1}{2}\Delta\alpha_1-\kappa^2_1+k^2-\frac{1}{4}(D\alpha_1)^2\Big)
(\hz\cdot a_1)\\
2(\hz\cdot D)I_3\hat R^H_1=2(\hz\cdot b_1)D\kappa_1\\
2(\hz\cdot D)\hat R^\Psi_1= -\Big(\frac{1}{2}\Delta\beta_1-\kappa^2_1+k^2-\frac{1}{4}(D\beta_1)^2\Big)
(\hz\cdot b_1)\\
2(\hz\cdot D)I_3\hat R^E_1=2(\hz\cdot a_1)D\kappa_1,
\end{array}
\right.
\end{equation}
\begin{equation}\label{equation_R_2_components}
\left\{\begin{array}{l}
2(\chz\cdot D)\check R^\Phi_2= -\Big(\frac{1}{2}\Delta\overline{\alpha}_2-\overline{\kappa}^2_2 +k^2-\frac{1}{4}(D\overline{\alpha}_2)^2\Big)
(\chz\cdot a_2)\\
2(\chz\cdot D)I_3\check R^H_2=2(\chz\cdot b_2)D\overline{\kappa}_2\\
2(\chz\cdot D)\check R^\Psi_2= -\Big(\frac{1}{2}\Delta\beta_2-\overline{\kappa}^2_2+k^2-\frac{1}{4}(D\beta_2)^2\Big)
(\chz\cdot b_2)\\
2(\chz\cdot D)I_3\check R^E_2=2(\chz\cdot a_2)D\overline{\kappa}_2.
\end{array}
\right.
\end{equation}
We take complex conjugates in (\ref{equation_R_2_components}) in order to get the information required.
\begin{equation}\label{equation_R_2_components_conjugate}
\left\{\begin{array}{l}
2(\hz\cdot D)\hat R^\Phi_2= \Big(\frac{1}{2}\Delta\alpha_2-\kappa^2_2 +k^2-\frac{1}{4}(D\alpha_2)^2\Big)(\hz\cdot \overline{a}_2)\\
2(\hz\cdot D)I_3\hat R^H_2=2(\hz\cdot \overline{b}_2)D\kappa_2\\
2(\hz\cdot D)\hat R^\Psi_2= \Big(\frac{1}{2}\Delta\beta_2-\kappa^2_2+k^2-\frac{1}{4}(D\beta_2)^2\Big)
(\hz\cdot \overline{b}_2)\\
2(\hz\cdot D)I_3\hat R^E_2=2(\hz\cdot \overline{a}_2)D\kappa_2.
\end{array}
\right.
\end{equation}
These equations are valid in $\mathbb{R}^3$. Multiplying the second and fourth equations of (\ref{equation_R_1_components}) and (\ref{equation_R_2_components_conjugate}) by $\hz\cdot$ one gets
\begin{equation}\label{enire_functions}
\left\{\begin{array}{l}
\hz\cdot \hat R^H_1=(\hz\cdot b_1)(\kappa_1-k)\\
\hz\cdot \hat R^E_1=(\hz\cdot a_1)(\kappa_1-k)\\
\hz\cdot \hat R^H_2=(\hz\cdot \overline{b}_2)(\kappa_2-k)\\
\hz\cdot \hat R^E_2=(\hz\cdot \overline{a}_2)(\kappa_2-k).
\end{array}
\right.
\end{equation}
To see the last identities we used the decay of $\hat R_j$ for $j=1,2$ at infinity, and the fact that $\hz\cdot \hat R^H_1-(\hz\cdot b_1)\kappa_1$, $\hz\cdot \hat R^E_1-(\hz\cdot a_1)\kappa_1$, $\hz\cdot \hat R^H_2-(\hz\cdot \overline{b}_2)\kappa_2$ and $\hz\cdot \hat R^E_2-(\hz\cdot \overline{a}_2)\kappa_2$ are entire functions in the variable $\hz\cdot x$.

Let us denote
\begin{eqnarray*}
(T1)&=&-\omega(\hz\cdot b_1)\int_O e^{i\xi\cdot x}
2(\hz\cdot D)\hat R^\Psi_2\left[\Big(\frac{\mu_1}{\mu_2}\Big)^{1/2}-1\right]\,dx\\
(T2)&=&-\omega(\hz\cdot a_1)\int_O e^{i\xi\cdot x}
2(\hz\cdot D)\hat R^\Phi_2 \left[\Big(\frac{\gamma_1}{\gamma_2}\Big)^{1/2}-1\right]\,dx\\
(T3)&=&-\omega(\hz\cdot\overline{b}_2)\int_O e^{i\xi\cdot x}
2(\hz\cdot D)\hat R^\Psi_1\left[\Big(\frac{\mu_1}{\mu_2}\Big)^{1/2}-1\right]\,dx\\
(T4)&=&-\omega(\hz\cdot\overline{a}_2)\int_O e^{i\xi\cdot x}
2(\hz\cdot D)\hat R^\Phi_1 \left[\Big(\frac{\gamma_1}{\gamma_2}\Big)^{1/2}-1\right]\,dx\\
\end{eqnarray*}
which are terms from (\ref{limit_computation_2}) and (\ref{limit_computation_3}). Then from the first and third equations in (\ref{equation_R_1_components}) and (\ref{equation_R_2_components_conjugate}) one gets
\begin{gather*}
(T1)+(T2)+(T3)+(T4)=\\
-\omega(\hz\cdot b_1)(\hz\cdot\overline{b}_2)\int_O e^{i\xi\cdot x} \Big[\frac{1}{2}\Delta(\beta_2-\beta_1)-\kappa^2_2+\kappa^2_1
-\frac{1}{4}(D\beta_2)^2+\frac{1}{4}(D\beta_1)^2\Big] \left[\Big(\frac{\mu_1}{\mu_2}\Big)^{1/2}-1\right]\,dx\\
-\omega(\hz\cdot a_1)(\hz\cdot\overline{a}_2)\int_O e^{i\xi\cdot x} \Big[\frac{1}{2}\Delta(\alpha_2-\alpha_1)-\kappa^2_2+\kappa^2_1
-\frac{1}{4}(D\alpha_2)^2+\frac{1}{4}(D\alpha_1)^2\Big] \left[\Big(\frac{\gamma_1}{\gamma_2}\Big)^{1/2}-1\right]\,dx.
\end{gather*}
Noting that
\begin{gather*}
\frac{1}{2}\Delta(\beta_2-\beta_1)=
D\Big(\frac{\mu_2}{\mu_1}\Big)^{1/2}\cdot D\Big(\frac{\mu_1}{\mu_2}\Big)^{1/2}
-\Big(\frac{\mu_2}{\mu_1}\Big)^{1/2}\Delta\Big(\frac{\mu_1}{\mu_2}\Big)^{1/2}\\
\frac{1}{2}\Delta(\alpha_2-\alpha_1)=
D\Big(\frac{\gamma_2}{\gamma_1}\Big)^{1/2}\cdot D\Big(\frac{\gamma_1}{\gamma_2}\Big)^{1/2}
-\Big(\frac{\gamma_2}{\gamma_1}\Big)^{1/2}\Delta\Big(\frac{\gamma_1}{\gamma_2}\Big)^{1/2},
\end{gather*}
one can write
\begin{gather*}
(T1)+(T2)+(T3)+(T4)=(T5)+(T6)\\
-\omega(\hz\cdot b_1)(\hz\cdot\overline{b}_2)\int_O e^{i\xi\cdot x} \Big[\Big(\frac{\mu_1}{\mu_2}\Big)^{1/2} D\Big(\frac{\mu_2}{\mu_1}\Big)^{1/2}\cdot D\Big(\frac{\mu_1}{\mu_2}\Big)^{1/2}
-\Delta\Big(\frac{\mu_1}{\mu_2}\Big)^{1/2} -\frac{1}{2}\Delta(\beta_2-\beta_1)\Big]\,dx\\
-\omega(\hz\cdot a_1)(\hz\cdot\overline{a}_2)\int_O e^{i\xi\cdot x} \Big[\Big(\frac{\gamma_1}{\gamma_2}\Big)^{1/2} D\Big(\frac{\gamma_2}{\gamma_1}\Big)^{1/2}\cdot D\Big(\frac{\gamma_1}{\gamma_2}\Big)^{1/2}
-\Delta\Big(\frac{\gamma_1}{\gamma_2}\Big)^{1/2} -\frac{1}{2}\Delta(\alpha_2-\alpha_1)\Big]\,dx,
\end{gather*}
where
\begin{gather*}
(T5)=-\omega(\hz\cdot b_1)(\hz\cdot\overline{b}_2)\int_O e^{i\xi\cdot x} \Big[-\kappa^2_2+\kappa^2_1-\frac{1}{4}(D\beta_2)^2+\frac{1}{4}(D\beta_1)^2\Big] \left[\Big(\frac{\mu_1}{\mu_2}\Big)^{1/2}-1\right]\,dx\\
(T6)=-\omega(\hz\cdot a_1)(\hz\cdot\overline{a}_2)\int_O e^{i\xi\cdot x} \Big[-\kappa^2_2+\kappa^2_1-\frac{1}{4}(D\alpha_2)^2+\frac{1}{4}(D\alpha_1)^2\Big] \left[\Big(\frac{\gamma_1}{\gamma_2}\Big)^{1/2}-1\right]\,dx.
\end{gather*}
Note that
\begin{gather*}
\Big(\frac{\mu_1}{\mu_2}\Big)^{1/2} D\Big(\frac{\mu_2}{\mu_1}\Big)^{1/2}\cdot D\Big(\frac{\mu_1}{\mu_2}\Big)^{1/2}=\frac{1}{2}D\beta_2\cdot D\Big(\frac{\mu_1}{\mu_2}\Big)^{1/2}-\frac{1}{2}D\beta_1\cdot D\Big(\frac{\mu_1}{\mu_2}\Big)^{1/2}\\
\Big(\frac{\gamma_1}{\gamma_2}\Big)^{1/2} D\Big(\frac{\gamma_2}{\gamma_1}\Big)^{1/2}\cdot D\Big(\frac{\gamma_1}{\gamma_2}\Big)^{1/2}=\frac{1}{2}D\alpha_2\cdot D\Big(\frac{\gamma_1}{\gamma_2}\Big)^{1/2}-\frac{1}{2}D\alpha_1\cdot D\Big(\frac{\gamma_1}{\gamma_2}\Big)^{1/2}.
\end{gather*}
Denoting
\begin{gather*}
(T7)=\omega(\hz\cdot b_1)(\hz\cdot\overline{b}_2)\int_O e^{i\xi\cdot x} \frac{1}{2}\Delta(\beta_2-\beta_1)\,dx\\
(T8)=\omega(\hz\cdot a_1)(\hz\cdot\overline{a}_2)\int_O e^{i\xi\cdot x}  \frac{1}{2}\Delta(\alpha_2-\alpha_1)\,dx
\end{gather*}
we get
\begin{gather*}
(T1)+(T2)+(T3)+(T4)=(T5)+(T6)+(T7)+(T8)\\
-\omega(\hz\cdot b_1)(\hz\cdot\overline{b}_2)\int_O e^{i\xi\cdot x} \Big[\frac{1}{2}D\beta_2\cdot D\Big(\frac{\mu_1}{\mu_2}\Big)^{1/2}-\frac{1}{2}D\beta_1\cdot D\Big(\frac{\mu_1}{\mu_2}\Big)^{1/2}
\Big]\,dx\\
-\omega(\hz\cdot a_1)(\hz\cdot\overline{a}_2)\int_O e^{i\xi\cdot x} \Big[\frac{1}{2}D\alpha_2\cdot D\Big(\frac{\gamma_1}{\gamma_2}\Big)^{1/2}-\frac{1}{2}D\alpha_1\cdot D\Big(\frac{\gamma_1}{\gamma_2}\Big)^{1/2}
\Big]\,dx\\
-\omega(\hz\cdot b_1)(\hz\cdot\overline{b}_2)\int_O e^{i\xi\cdot x}
(-\Delta)\Big(\frac{\mu_1}{\mu_2}\Big)^{1/2}\,dx\\
-\omega(\hz\cdot a_1)(\hz\cdot\overline{a}_2)\int_O e^{i\xi\cdot x}
(-\Delta)\Big(\frac{\gamma_1}{\gamma_2}\Big)^{1/2}\,dx.
\end{gather*}
We consider the next terms from (\ref{limit_computation_1}) and (\ref{limit_computation_2})
\begin{gather*}
(T9)=\omega(\hz\cdot\overline{b}_2)(\hz\cdot b_1)\int_O e^{i\xi\cdot x}
(-\Delta)\Big(\frac{\mu_1}{\mu_2}\Big)^{1/2}\,dx\\
(T10)=\omega(\hz\cdot\overline{a}_2)(\hz\cdot a_1)\int_O e^{i\xi\cdot x}
(-\Delta)\Big(\frac{\gamma_1}{\gamma_2}\Big)^{1/2}\,dx\\
(T11)=\omega\int_O e^{i\xi\cdot x}\left[(\hz\cdot \overline{a}_2)(\hz\cdot a_1) D\alpha_2\cdot D\Big(\frac{\gamma_1}{\gamma_2}\Big)^{1/2}+
(\hz\cdot \overline{b}_2)(\hz\cdot b_1)
D\beta_2\cdot D\Big(\frac{\mu_1}{\mu_2}\Big)^{1/2}\right]\,dx,
\end{gather*}
hence
\begin{gather*}
(T1)+(T2)+(T3)+(T4)+(T9)+(T10)+(T11)=(T5)+(T6)+(T7)+(T8)\\
+\omega(\hz\cdot b_1)(\hz\cdot\overline{b}_2)\int_O e^{i\xi\cdot x} \Big[\frac{1}{2}D\beta_2\cdot D\Big(\frac{\mu_1}{\mu_2}\Big)^{1/2}+\frac{1}{2}D\beta_1\cdot D\Big(\frac{\mu_1}{\mu_2}\Big)^{1/2}
\Big]\,dx\\
+\omega(\hz\cdot a_1)(\hz\cdot\overline{a}_2)\int_O e^{i\xi\cdot x} \Big[\frac{1}{2}D\alpha_2\cdot D\Big(\frac{\gamma_1}{\gamma_2}\Big)^{1/2}+\frac{1}{2}D\alpha_1\cdot D\Big(\frac{\gamma_1}{\gamma_2}\Big)^{1/2}
\Big]\,dx.
\end{gather*}
We note that
\begin{equation*}
D\Big(\frac{\mu_1}{\mu_2}\Big)^{1/2}=
\frac{1}{2}\Big(\frac{\mu_1}{\mu_2}\Big)^{1/2}(D\beta_1-D\beta_2),\qquad
D\Big(\frac{\gamma_1}{\gamma_2}\Big)^{1/2}=
\frac{1}{2}\Big(\frac{\gamma_1}{\gamma_2}\Big)^{1/2}(D\alpha_1-D\alpha_2)
\end{equation*}
hence
\begin{gather}
(T1)+(T2)+(T3)+(T4)+(T9)+(T10)+(T11)=(T5)+(T6)+(T7)+(T8) \notag \\
+\omega(\hz\cdot b_1)(\hz\cdot\overline{b}_2)\int_O e^{i\xi\cdot x} \frac{1}{4}\Big(\frac{\mu_1}{\mu_2}\Big)^{1/2}[(D\beta_1)^2-(D\beta_2)^2]\,dx \notag \\
+\omega(\hz\cdot a_1)(\hz\cdot\overline{a}_2)\int_O e^{i\xi\cdot x} \frac{1}{4}\Big(\frac{\gamma_1}{\gamma_2}\Big)^{1/2}[(D\alpha_1)^2-(D\alpha_2)^2]\,dx. \label{finalterms1}
\end{gather}
Considering the terms
\begin{gather*}
(T12)=\int_O e^{i\xi\cdot x}
\left[(\hz\cdot\hat R^E_2)(\hz\cdot a_1) (\kappa_2\tga+\kappa_1\tmu)+
(\hz\cdot\hat R^H_2)(\hz\cdot b_1) (\kappa_2\tmu+\kappa_1\tga)\right]\,dx\\
(T13)=\int_O e^{i\xi\cdot x}
\left[(\kappa_2\tga+\kappa_1\tmu)(\hz\cdot\overline{b}_2)(\hz\cdot\hat R^H_1)
+(\kappa_2\tmu+\kappa_1\tga)(\hz\cdot\overline{a}_2)(\hz\cdot\hat R^E_1)\right]\,dx\\
(T14)=k((\hz\cdot\overline{b}_2)(\hz\cdot b_1)+ (\hz\cdot\overline{a}_2)(\hz\cdot a_1))\int_O e^{i\xi\cdot x} (\kappa_1+\kappa_2)(\tmu+\tga)\,dx\\
(T15)=-\int_O e^{i\xi\cdot x}\kappa_2\left[(\hz\cdot \overline{a}_2)(\hz\cdot a_1) (\kappa_2\tga+\kappa_1\tmu)+(\hz\cdot \overline{b}_2)(\hz\cdot b_1) (\kappa_2\tmu+\kappa_1\tga)\right]\,dx
\end{gather*}
from (\ref{limit_computation_1}), (\ref{limit_computation_2}) and (\ref{limit_computation_3}); and using (\ref{enire_functions}) we have that
\begin{gather*}
(T12)+(T13)=\\
=\int_O e^{i\xi\cdot x}
\left[(\hz\cdot \overline{a}_2)(\hz\cdot a_1)(\kappa_2-k)(\kappa_2\tga+\kappa_1\tmu)+
(\hz\cdot \overline{b}_2)(\hz\cdot b_1)(\kappa_2-k)(\kappa_2\tmu+\kappa_1\tga)\right]\,dx\\
+\int_O e^{i\xi\cdot x}
\left[(\kappa_2\tga+\kappa_1\tmu)(\hz\cdot\overline{b}_2)(\hz\cdot b_1)(\kappa_1-k)
+(\kappa_2\tmu+\kappa_1\tga)(\hz\cdot\overline{a}_2)(\hz\cdot a_1)(\kappa_1-k)\right]\,dx.
\end{gather*}
Therefore
\begin{gather}
(T12)+(T13)+(T14)+(T15)= \notag \\
=(\hz\cdot \overline{a}_2)(\hz\cdot a_1)\int_O e^{i\xi\cdot x}
\kappa_1(\kappa_2\tmu+\kappa_1\tga)\,dx
+(\hz\cdot\overline{b}_2)(\hz\cdot b_1)\int_O e^{i\xi\cdot x}
\kappa_1(\kappa_2\tga+\kappa_1\tmu)\,dx. \label{finalterms2}
\end{gather}

Now we make the choice $b_1=\overline{b}_2=\chz$ and $a_1=\overline{a}_2=\hz$. The equations (\ref{equation_R_1_components}) and (\ref{equation_R_2_components_conjugate}) imply that $\hat R^E_j=0$ for $j=1,2$. Considering all the terms from the technical lemmas in Section \ref{sec:technical}, we see that the only surviving terms in \eqref{limit_orthogonal_relation_1}--\eqref{limit_crossoed_terms} are those which appear on the left hand sides of \eqref{finalterms1} and \eqref{finalterms2}. Thus these computations give
\begin{eqnarray*}
\int_O e^{i\xi\cdot x} \left[\frac{1}{2}\Delta(\beta_2-\beta_1)+ \frac{1}{4}\Big(\frac{\mu_1}{\mu_2}\Big)^{1/2}[(D\beta_1)^2-(D\beta_2)^2]
+\frac{\kappa_1}{\omega}(\kappa_2\tga+\kappa_1\tmu)\right.\\
\left.+(\kappa^2_2-\kappa^2_1+\frac{1}{4}(D\beta_2)^2-\frac{1}{4}(D\beta_1)^2)
\Big[\Big(\frac{\mu_1}{\mu_2}\Big)^{1/2}-1\Big]\right]\,dx=0.
\end{eqnarray*}
If we choose $a_1=\overline{a}_2=\chz$ and $b_1=\overline{b}_2=\hz$, we get that $\hat R^H_j=0$ for $j=1,2$ for the same reason as above. We consider again all the terms from the lemmas in Section \ref{sec:technical} and we get
\begin{eqnarray*}
\int_O e^{i\xi\cdot x} \left[\frac{1}{2}\Delta(\alpha_2-\alpha_1)+ \frac{1}{4}\Big(\frac{\gamma_1}{\gamma_2}\Big)^{1/2}[(D\alpha_1)^2-(D\alpha_2)^2]
+\frac{\kappa_1}{\omega}(\kappa_2\tmu+\kappa_1\tga)\right.\\
\left.+(\kappa^2_2-\kappa^2_1+\frac{1}{4}(D\alpha_2)^2-\frac{1}{4}(D\alpha_1)^2)
\Big[\Big(\frac{\gamma_1}{\gamma_2}\Big)^{1/2}-1\Big]\right]\,dx=0.
\end{eqnarray*}
Finally, by a direct computation we see that 
\begin{eqnarray*}
 & \frac{\kappa_1}{\omega}(\kappa_2 \tilde{\mu} + \kappa_1 \tilde{\gamma}) + (\kappa_2^2-\kappa_1^2) \big( \frac{\gamma_1}{\gamma_2} \big)^{1/2} = 0, & \\
 & \frac{\kappa_1}{\omega}(\kappa_2 \tilde{\gamma} + \kappa_1 \tilde{\mu}) + (\kappa_2^2-\kappa_1^2) \big( \frac{\mu_1}{\mu_2} \big)^{1/2} = 0, & 
\end{eqnarray*}
hence the result is obtained.

The preceding arguments are valid for all $\xi \in \R^3$ with $|\xi'| > 0$, since this was the assumption on $\xi$ in the construction of reflected CGO solutions. However, the continuity of the Fourier transform proves the result for all $\xi$.
\end{proof}
\medskip

The above proposition shows that the coefficients satisfy the following equations in $O$: 
\begin{eqnarray*}
 & -\frac{1}{2} \Delta (\alpha_1-\alpha_2) + \frac{1}{4}D(\alpha_1+\alpha_2) \cdot D(\alpha_1-\alpha_2) + \omega^2 (\mu_1 \gamma_1 - \mu_2 \gamma_2) = 0, & \\
 & -\frac{1}{2} \Delta (\beta_1-\beta_2) + \frac{1}{4}D(\beta_1+\beta_2) \cdot D(\beta_1-\beta_2) + \omega^2 (\mu_1 \gamma_1 - \mu_2 \gamma_2) = 0. & 
\end{eqnarray*}
Let $u = (\gamma_1/\gamma_2)^{1/2}$ and $v = (\mu_1/\mu_2)^{1/2}$. The equations become 
\begin{eqnarray*}
 & - \Delta (\log\,u) + (\gamma_1 \gamma_2)^{-1/2} D(\gamma_1 \gamma_2)^{1/2} \cdot D(\log\,u) + \omega^2 (\mu_1 \gamma_1 - \mu_2 \gamma_2) = 0, & \\
 & - \Delta (\log\,v) + (\mu_1 \mu_2)^{-1/2} D(\mu_1 \mu_2)^{1/2} \cdot D(\log\,v) + \omega^2 (\mu_1 \gamma_1 - \mu_2 \gamma_2) = 0. & 
\end{eqnarray*}
Multiplying the equations by $(\gamma_1 \gamma_2)^{1/2}$ and $(\mu_1 \mu_2)^{1/2}$ we obtain 
\begin{eqnarray*}
 & D \cdot ((\gamma_1 \gamma_2)^{1/2} D(\log\,u)) + \omega^2 (\mu_1 \gamma_1 - \mu_2 \gamma_2)(\gamma_1 \gamma_2)^{1/2} = 0, & \\
 & D \cdot ((\mu_1 \mu_2)^{1/2} D(\log\,v)) + \omega^2 (\mu_1 \gamma_1 - \mu_2 \gamma_2)(\mu_1 \mu_2)^{1/2} = 0. & 
\end{eqnarray*}
Finally, these may be written as 
\begin{eqnarray*}
 & D \cdot (\gamma_2 Du) + \omega^2 \mu_2 \gamma_2^2 (u^2 v^2 - 1)u = 0, & \\
 & D \cdot (\mu_2 Dv) + \omega^2 \mu_2^2 \gamma_2(u^2 v^2 - 1)v = 0 & 
\end{eqnarray*}
in $O$. This is a semilinear elliptic system for $u, v \in C^4(\overline{O})$. The assumptions in Theorem 1 imply that $u|_{\partial O} = v|_{\partial O} = 1$ and $\partial_{\nu} u|_{\partial O} = \partial_{\nu} v|_{\partial O} = 0$. The following formulation of the unique continuation property will therefore imply that $\gamma_1 \equiv \gamma_2$ and $\mu_1 \equiv \mu_2$ in $O$, thus proving Theorem 1.

\begin{lemma}
Let $a$, $b$ be non-vanishing complex functions in $C^2(\closure{O})$ with positive real part, and let $p, q$ be complex functions in $L^{\infty}(O)$. Suppose that $u, v \in C^2(\closure{O})$ satisfy in $O$ 
\begin{eqnarray*}
 & D \cdot (a Du) + p (u^2 v^2 - 1)u = 0, & \\
 & D \cdot (b Dv) + q (u^2 v^2 - 1)v = 0, & 
\end{eqnarray*}
and also $u|_{\partial O} = v|_{\partial O} = 1$ and $\partial_{\nu} u|_{\partial O} = \partial_{\nu} v|_{\partial O} = 0$. Then $u \equiv v \equiv 1$ in $O$.
\end{lemma}
\begin{proof}
Let $u = 1 + a^{-1/2} z$ and $v = 1 + b^{-1/2} w$. Using the identities 
\begin{eqnarray*}
 & D \cdot (a D(a^{-1/2} z)) = a^{1/2}(-\Delta z + \frac{\Delta(a^{1/2})}{a^{1/2}} z), & \\
 & D \cdot (b D(b^{-1/2} w)) = b^{1/2}(-\Delta w + \frac{\Delta(b^{1/2})}{b^{1/2}} w), & \\
 & u^2 v^2 - 1 = (ab)^{-1/2} (zw + b^{1/2} z + a^{1/2} w)(uv+1), & 
\end{eqnarray*}
we obtain the equations 
\begin{eqnarray*}
 & -\Delta z + \frac{\Delta(a^{1/2})}{a^{1/2}} z + \tilde{p} (zw + b^{1/2} z + a^{1/2} w) = 0, & \\
 & -\Delta w + \frac{\Delta(b^{1/2})}{b^{1/2}} w + \tilde{q} (zw + b^{1/2} z + a^{1/2} w) = 0. & 
\end{eqnarray*}
Here $\tilde{p}, \tilde{q} \in L^{\infty}(O)$ are given by 
\begin{eqnarray*}
 & \tilde{p} = a^{-1/2} p (ab)^{-1/2} (uv+1)u, & \\
 & \tilde{q} = b^{-1/2} q (ab)^{-1/2} (uv+1)v. & 
\end{eqnarray*}
Thus, there is $C > 0$ such that the following differential inequalities for $z, w \in C^2(\closure{O})$ are valid almost everywhere in $O$:
\begin{eqnarray*}
 & \abs{\Delta z(x)} \leq C (\abs{z(x)} + \abs{w(x)}), & \\
 & \abs{\Delta w(x)} \leq C (\abs{z(x)} + \abs{w(x)}). & 
\end{eqnarray*}
We also have $z|_{\partial O} = w|_{\partial O} = 0$ and $\partial_{\nu} z|_{\partial O} = \partial_{\nu} w|_{\partial O} = 0$. The unique continuation property holds in this setting, which can be seen by applying a scalar Carleman estimate to both $z$ and $w$ (for details see \cite{KSaU}). Thus $z$ and $w$ vanish identically in $O$.
\end{proof}

\medskip

\section{Proof of Theorem 2.}

Here we prove an analog of Theorem 1 for the case  where the inaccessible part is part of a sphere. Following the idea in \cite{I}, this reduces to the hyperplane case by using the Kelvin transform.

We first need to check how the Maxwell equations behave under the Kelvin transform. For this it is convenient to consider $E$ and $H$ as complex valued $1$-forms in $\Omega$ (that is, a vector field $(E_1,E_2,E_3)$ is identified with the $1$-form $E_1 \,dx^1 + E_2 \,dx^2 + E_3 \,dx^3$). We equip $\Omega$ with the Euclidean metric $e$, and introduce the Hodge star operator $*_e$ which maps a $1$-form $X_1 \,dx^1 + X_2 \,dx^2 + X_3 \,dx^3$ to the $2$-form $X_1 \,dx^2 \wedge dx^3 + X_2 \,dx^3 \wedge dx^1 + X_3 \,dx^1 \wedge dx^2$. If $d$ is the exterior derivative, the Maxwell equations in form notation are 
\begin{equation*}
\left\{ \begin{array}{rl}
d E &= i \omega \mu *_e H \\
d H &= -i \omega \gamma *_e E
\end{array} \right. \quad \text{in } \Omega.
\end{equation*}

Assume that $0$ is not in $\Omega$. The Kelvin transform is given by 
\begin{equation*}
\tilde{x}(x) = \abs{x}^{-2} x, \quad F(\tilde{x}) = x(\tilde{x}) = \abs{\tilde{x}}^{-2} \tilde{x}.
\end{equation*}
If $\tilde{\Omega} = F^{-1}(\Omega)$, then $F$ is a conformal transformation from $(\tilde{\Omega}, e)$ onto $(\Omega,e)$: 
\begin{equation*}
F^* e = \abs{\tilde{x}}^{-4} e.
\end{equation*}
Let $\tilde{E} = F^* E$, $\tilde{H} = F^* H$, $\tilde{\mu} = F^* \mu$, and $\tilde{\gamma} = F^* \gamma$. The following is the transformation law for the Maxwell equations under the Kelvin transform.

\begin{lemma}
One has 
\begin{equation*}
\left\{ \begin{array}{rl}
d E &= i \omega \mu *_e H \\
d H &= -i \omega \gamma *_e E
\end{array} \right. \ \ \text{in } \Omega \quad \Leftrightarrow \quad \left\{ \begin{array}{rl}
d \tilde{E} &= i \omega \tilde{\mu} \abs{\tilde{x}}^{-2} *_e \tilde{H} \\
d \tilde{H} &= -i \omega \tilde{\gamma} \abs{\tilde{x}}^{-2} *_e \tilde{E}
\end{array} \right. \ \ \text{in } \tilde{\Omega}.
\end{equation*}
\end{lemma}
\begin{proof}
We use the following identities valid for $k$-forms $\eta$ in a $3$-manifold:
\begin{equation*}
d F^* \eta = F^* d\eta, \quad F^* (*_e \eta) = *_{F^* e} F^* \eta, \quad *_{ce} \eta = c^{3/2-k} *_e \eta.
\end{equation*}
Here $c$ is any positive smooth function. If $dE = i\omega \mu *_e H$, then 
\begin{align*}
d\tilde{E} &= d F^* E = F^* dE = F^* (i\omega \mu *_e H) = i\omega \tilde{\mu} *_{F^*e} \tilde{H} = i\omega \tilde{\mu} *_{\abs{\tilde{x}}^{-4} e} \tilde{H} \\
 &= i\omega \tilde{\mu} \abs{\tilde{x}}^{-2} *_e \tilde{H}.
\end{align*}
The other direction and the other equation are handled analogously.
\end{proof}

\medskip
Next we check how the impedance map transforms. Note that in the form notation, the boundary condition in the Maxwell equations corresponds to fixing the tangential $1$-form $*_e (\nu \wedge H)$ on $\partial \Omega$, where $\nu = \nu_1 \,dx^1 + \nu_2 \,dx^2 + \nu_3 \,dx^3$ is the outer unit normal written as a $1$-form.

\begin{lemma} \label{lemma:impedance_map_kelvin}
Let $\Lambda$ be the impedance map in $\Omega$ with coefficients $(\mu,\gamma)$, and let $\tilde{\Lambda}$ be the impedance map in $\tilde{\Omega}$ with coefficients $(\tilde{\mu} \abs{\tilde{x}}^{-2}, \tilde{\gamma} \abs{\tilde{x}}^{-2})$. If $\tilde{T}$ is a tangential field on $\partial \tilde{\Omega}$, then 
\begin{equation*}
\tilde{\Lambda}(\tilde{T}) = F^* \Lambda( (F^{-1})^* \tilde{T}).
\end{equation*}
\end{lemma}
\begin{proof}
We take $\rho$ to be a boundary defining function for $\Omega$, that is, $\rho$ is a smooth function $\R^3 \to \R$ and $\Omega = \{ \rho > 0 \}$, $\partial \Omega = \{ \rho = 0 \}$. Then $\tilde{\rho} = F^* \rho$ is a boundary defining function for $\tilde{\Omega}$. The unit normal is the $1$-form 
\begin{equation*}
\tilde{\nu} = -\frac{d\tilde{\rho}}{\abs{d\tilde{\rho}}_e}  = -\frac{F^* d \rho}{\abs{F^* d \rho}_{F^* (F^{-1})^* e}} = - F^* \frac{d\rho}{\abs{d\rho}_{(F^{-1})^* e}} = -F^* \left( \lvert x \rvert^{-2} \frac{d\rho}{\lvert d\rho \rvert_{e}} \right) = \abs{\tilde{x}}^2 F^* \nu,
\end{equation*}
using that $(F^{-1})^* e = \abs{x}^{-4} e$ and $\nu = -d\rho/\abs{d\rho}_e$.

Let $(\tilde{E},\tilde{H})$ be a solution to the Maxwell system in $\tilde{\Omega}$ with $*_e (\tilde{\nu} \wedge \tilde{H})|_{\partial \tilde{\Omega}} = \tilde{T}$. Then on $\partial \tilde{\Omega}$ one has 
\begin{equation*}
*_e (\tilde{\nu} \wedge \tilde{E}) = \abs{\tilde{x}}^2 *_e F^* (\nu \wedge E) = \abs{\tilde{x}}^2 F^* ( *_{(F^{-1})^* e} (\nu \wedge E) ) = F^* ( *_e (\nu \wedge E)).
\end{equation*}
Similarly $*_e (\tilde{\nu} \wedge \tilde{H}) = F^* ( *_e (\nu \wedge H))$, and therefore 
\begin{equation*}
\tilde{\Lambda}(\tilde{T}) = *_e (\tilde{\nu} \wedge \tilde{E}) = F^* \Lambda (*_e (\nu \wedge H)) = F^* \Lambda ( (F^{-1})^* \tilde{T} ).
\end{equation*}
\end{proof}

\medskip
We now assume that $(\mu_j, \gamma_j)$ are two sets of coefficients satisfying the assumptions of Theorem 2, and let \(\Lambda _1\) and \(\Lambda _2\) be the corresponding impedance maps. Let $B_0$ be an open ball with $\Omega \subset B_0$, and suppose that $\Gamma \subset \partial \Omega$ is such that $\Lambda_1 T|_{\Gamma} = \Lambda_2 T|_{\Gamma}$ for tangential fields $T$ supported in $\Gamma$. Assume that $\Gamma_0 = \partial \Omega \setminus \Gamma$ satisfies $\Gamma_0 = \partial \Omega \cap \partial B_0$ and $\Gamma_0 \neq \partial B_0$, that is, the inaccessible part $\Gamma_0$ is part of a sphere.
\smallskip

We wish to show that $\mu_1 = \mu_2$ and $\gamma_1 = \gamma_2$ in $\Omega$. To this end, choose coordinates so that $B_0 = B(x_0,1/2)$ where $x_0 = (0,0,1/2)$, and assume that the origin is not in $\closure{\Omega}$. Let $\tilde{\Omega} = F^{-1}(\Omega)$ be the image of $\Omega$ under the Kelvin transform. From Lemma \ref{lemma:impedance_map_kelvin} we obtain that 
\begin{equation*}
\tilde{\Lambda}_1(\tilde{T})|_{\tilde{\Gamma}} = \tilde{\Lambda}_2(\tilde{T})|_{\tilde{\Gamma}}
\end{equation*}
for tangential fields supported in $\tilde{\Gamma} = F^{-1}(\Gamma)$. Since the Kelvin transform maps $\Gamma_0$ onto a subset of $\{x_3 = 1\}$, we are in a situation where the inaccessible part of the boundary is part of a hyperplane. Also the other conditions in Theorem 1 are satisfied; in particular, \eqref{BoundaryExtensionAssumption} follows from its analog for $\partial B_0$ which reads  
\begin{equation*}
\left\{ \begin{array}{l} \text{there exist $C^4$ extensions of $\gamma_j$ and $\mu_j$ into $\R^3$ which are } \\[1pt] \text{preserved under the map $x \mapsto F \circ R \circ F^{-1}(x)$}, \end{array} \right.
\end{equation*}
where $F$ is the Kelvin transform and $R$ is the reflection $(\tilde{x}_1,\tilde{x}_2,\tilde{x}_3) \mapsto (\tilde{x}_1,\tilde{x}_2,2-\tilde{x}_3)$. Consequently, Theorem 1 implies that 
\begin{equation*}
\tilde{\mu}_1 \abs{\tilde{x}}^{-2} = \tilde{\mu}_2 \abs{\tilde{x}}^{-2}, \quad \tilde{\gamma}_1 \abs{\tilde{x}}^{-2} = \tilde{\gamma}_2 \abs{\tilde{x}}^{-2} \quad \text{in } \tilde{\Omega}.
\end{equation*}
It follows that $\mu_1 = \mu_2$ and $\gamma_1 = \gamma_2$ in $\Omega$, thus proving Theorem 2.



\begin{thebibliography}{XXXX}

\bibitem[BU]{B-U}
Bukhgeim, A.~L., Uhlmann, G., {Recovering a potential from partial
  Cauchy data}. \textit{Comm. PDE} \textbf{27} (2002), p.~653--668.

\bibitem[C]{C} Calder\'on, A.~P., On an inverse boundary value problem. Seminar on Numerical Analysis and its Applications to Continuum Physics (R{\'i}o de Janeiro, 1980), pp. 65--73, \textit{Soc. Brasil. Mat.}, R{\'i}o de Janeiro, 1980.

\bibitem[CP]{C-P}
Colton, D., P\"aiv\"arinta, L., The uniqueness of a solution to an inverse scattering problem for electromagnetic waves. \textit{Arch. Rational Mech. Anal.} \textbf{119} (1992), no. 1, 59--70.

\bibitem[I]{I} Isakov, V., On uniqueness in the inverse conductivity problem with local data. \textit{Inverse Probl. Imaging} \textbf{1} (2007), no. 1, p. 95-105.

\bibitem[JM]{J-M} Joshi, M., McDowall, S.~R., Total determination of material parameters from electromagnetic boundary information. \textit{Pacific J. Math.} \textbf{193} (2000), 107--129.

\bibitem[KSaU]{KSaU} Kenig, C.~E., Salo, M., Uhlmann, G., Inverse problems for the anisotropic Maxwell equations. Preprint (2009).

\bibitem[KSU]{K-S-U} Kenig, C.~E., Sj\"ostrand J., Uhlmann, G., The Calder\'on problem with partial data. \textit{Ann. of Math.} \textbf{165} (2007), p.~567--591.

\bibitem[McD1]{McD1}
McDowall, S.~R., Boundary determination of material parameters from electromagnetic boundary information. \textit{Inverse Problems} \textbf{13} (1997), 153--163.

\bibitem[McD2]{McD2}
McDowall, S.~R., An electromagnetic inverse problem in chiral media. \textit{Trans. Amer. Math. Soc.} \textbf{352} (2000), no. 7, 2993--3013.

\bibitem[MM]{M-M} Mitrea, D., Mitrea, M.,
Finite energy solutions of Maxwell's equations and constructive Hodge decompositions on nonsmooth Riemannian manifolds. \textit{J. Funct. Anal.} \textbf{190} (2002), no. 2, 339--417.

\bibitem[N]{N} Nachman, A.~I., Reconstructions from boundary measurements. \textit{Ann. of Math.} \textbf{128} (1988), no. 3, 531--576.

\bibitem[OPS1]{O-P-S1}
Ola, P., P\"aiv\"arinta, L., Somersalo, E., An inverse boundary value problem in electrodynamics. \textit{Duke Math. J.} \textbf{70} (1993), no. 3, 617--653.

\bibitem[OPS2]{O-P-S2}
Ola, P., P\"aiv\"arinta, L., Somersalo, E., Inverse problems for time harmonic electrodynamics. Inside out: inverse problems and applications, 169--191, \textit{Math. Sci. Res. Inst. Publ.}, \textbf{47}, Cambridge Univ. Press, Cambridge, 2003.
 
\bibitem[OS]{O-S} Ola, P., Somersalo, E., Electromagnetic inverse problems and generalized Sommerfeld potentials. \textit{SIAM J. Appl. Math.} \textbf{56} (1996), no. 4, 1129-1145.

\bibitem[S]{S} Salo, M., Semiclassical pseudodifferential calculus and the reconstruction of a magnetic field. \textit{Comm. PDE} \textbf{31} (2006), no. 11, p. 1639-1666.

\bibitem[ST]{SaT}
Salo, M., Tzou, L., Inverse problems with partial data for a Dirac system: a Carleman estimate approach. Preprint (2009).

\bibitem[S1]{S1}
Saranen, J., \"Uber das Verhalten der L\"osungen der Maxwellschen Randwertaufgabe in Gebieten mit Kegelspitzen. \textit{Math. Methods Appl. Sci.}  \textbf{2} (1980), no. 2, 235--250. 

\bibitem[S2]{S2} 
Saranen, J.,  \"Uber das Verhalten der L\"osungen der Maxwellschen Randwertaufgabe in einigen nichtglatten Gebieten. \textit{Ann. Acad. Sci. Fenn. Ser. A I Math.} \textbf{6} (1981), no. 1, 15--28.

\bibitem[Sa]{Sa}
Sarkola, E., A unified approach to direct and inverse scattering for acoustic and electromagnetic waves. 
\textit{Ann. Acad. Sci. Fenn. Math. Diss.}  \textbf{101} (1995).

\bibitem[SU]{S-U} Sylvester, J., Uhlmann, G., A global uniqueness theorem for an inverse boundary value problem. \textit{Ann. of Math.} \textbf{125} (1987), 153-169.

\end{thebibliography}
\end{document}